\documentclass[article,12pt,letterpaper]{amsart}

\usepackage{amsmath,amssymb,amsthm,graphicx,mathrsfs,url}
\usepackage[usenames,dvipsnames]{color}
\usepackage{hyperref}
\usepackage{amsxtra}
\usepackage{graphicx}
\usepackage{tikz}
\usepackage{xcolor}
\relax

\numberwithin{equation}{section}

\def\11{{\rm 1~\hspace{-1.4ex}l} }
\def\R{\mathbb R}

\def\Z{\mathbb Z}
\def\N{\mathbb N}

\def\T{\mathbb T}

\renewcommand{\H}{\mathbf{H}}

\newcommand{\e}{\varepsilon}

\newcommand{\eps}{{\epsilon}}

\renewcommand{\div}{\text{div}}

\theoremstyle{plain}

\newtheorem{thm}{Theorem}
\newtheorem{prop}{Proposition}[section]
\newtheorem{cor}[prop]{Corollary}

\newtheorem{definition}[prop]{Definition}

\newtheorem{remark}[prop]{Remark}

\theoremstyle{definition}

\newtheorem{example}[prop]{Example}
\newtheorem{rem}[prop]{Remark}
 
\numberwithin{equation}{section}

\definecolor{darkred}{rgb}{0.7,0.1,0.1}
\definecolor{darkblue}{rgb}{0,0,0.7}
\title[Exact observability of subelliptic PDEs]{Exact observability properties of subelliptic wave and Schr\"odinger equations}
 \author[C. Letrouit]{Cyril Letrouit}
\address[C. Letrouit]{Sorbonne Universit\'e, Universit\'e Paris-Diderot, CNRS, Inria, Laboratoire Jacques-Louis Lions, F-75005 Paris. \newline DMA, \'Ecole normale sup\'erieure, CNRS, PSL Research University, 75005 Paris.}
\email{letrouit@ljll.math.upmc.fr}

\begin{document}    

\maketitle

\begin{abstract}
In this survey paper, we report on recent works concerning exact observability (and, by duality, exact controllability) properties of subelliptic wave and Schr\"odinger-type equations. These results illustrate the slowdown of propagation in directions transverse to the horizontal distribution. The proofs combine sub-Riemannian geometry, semi-classical analysis, spectral theory and non-commutative harmonic analysis. 
\end{abstract}
\setcounter{tocdepth}{1}
\tableofcontents

\section{Introduction}
\subsection{Controllability and observability} The problem of (exact) controllability of PDEs, which has been intensively studied in the past decades, is the following: given a manifold $M$, a subset $\omega\subset M$, a time $T>0$ and an operator $A$ acting on functions on $M$, the study of exact controllability consists in determining whether, for any initial state $u_0$ and any final state $u_1$, there exists $f$ such that the solution of
\begin{equation} \label{e:control}
\partial_tu=Au+\mathbf{1}_{\omega}f, \qquad u_{|t=0}=u_0
\end{equation}
in $M$ is equal to $u_1$ at time $T$. Here, $\mathbf{1}_{\omega}$ is the characteristic function of $\omega$. In other words, exact controllability holds if it is possible, starting from any initial state, to reach any final state just acting on $\omega$ during a time $T$. The general answer depends on  the time $T$, the control set $\omega$, the operator $A$, and the functional spaces in which $u_0$, $u_1$ and $f$ live. This problem is relevant in many physical situations: typical examples are the control of the temperature of a room by a heater, or the acoustic insulation of a room just by acting on a small part of it.

By duality (Hilbert Uniqueness Method, see \cite{lions1988controlabilite}), the exact controllability property is equivalent to some inequality of the form
\begin{equation}\label{e:generobs}
\exists C_{T,\omega}>0, \ \forall u_0, \qquad \|u_0\|^2 \leq C_{T,\omega}\int_0^T \|\mathbf{1}_\omega u(t)\|^2dt, 
\end{equation}
where $u$ is the solution of the adjoint equation $(\partial_t+A^*)u=0$ with initial datum $u_0$ (here again, one should specify functional spaces). This is called an observability inequality.  In other words, controllability holds if and only if any solution of $(\partial_t+A^*)u=0$ can be detected from $\omega$, in a ``quantitative way'' which is measured by the constant $C_{T,\omega}$. This paper is devoted to the study of observability for equations of wave-type, Schr\"odinger-type or heat-type, i.e. we consider the equations
\begin{align}
(\partial^2_{tt}-\mathcal{L})u=0 \qquad \qquad&\text{(Wave-type)}, \label{e:waveop} \\
(i\partial_t-\mathcal{L})u=0 \qquad \qquad&\text{(Schr\"odinger-type)}, \label{e:schrodop} \\
(\partial_{t}-\mathcal{L})u=0 \qquad \qquad&\text{(Heat-type)} \label{e:heatop}
\end{align}
for various time-independent operators $\mathcal{L}$ on $M$.\footnote{The wave equation involves a $\partial_{tt}^2$ term, and thus does not enter, strictly speaking, the framework given by equation \ref{e:control}. However, it is possible to give a common framework for all three equations, at the cost of being a bit more abstract. See \cite[Section 2.3]{coron2007control} for a general introduction.} By duality, all the observability results presented here imply exact controllability results as explained above, but we won't state them for the sake of simplicity.

\subsection{Observability of classical PDEs} Let us present a first series of results, dating back to the 1990's, which concern the observability problem in case $M$ is a compact Riemannian manifold with a metric $g$ and with boundary $\partial M\neq\emptyset$, $\mathcal{L}=\Delta_g$ is the Laplace-Beltrami operator on $(M,g)$ and the equation is one of the three equations \eqref{e:waveop}, \eqref{e:schrodop} or \eqref{e:heatop}, with Dirichlet boundary conditions $u|_{\partial M}=0$. We deal with these three problems in this order, following the chronology of the results.

Throughout this section, $(M,g)$ is a fixed manifold with boundary $\partial M\neq\emptyset$ and $\mathcal{L}=\Delta_g$. In this section, the notation $dx$ stands for the associated Riemannian volume $dx=d{\rm vol}_g(x)$.

\subsubsection{Observability of the Riemannian wave equation} Let us start with the wave equation \eqref{e:waveop} with initial data $(u_{t=0},\partial_tu_{|t=0})=(u_0,u_1)\in H^1(M)\times L^2(M)$ and Dirichlet boundary conditions. The energy of a solution, which is conserved along the flow, is 
\begin{equation*}
E(u(t))=\|(-\Delta_g)^{\frac12}u(t,\cdot)\|^2_{L^2(M)}+\|\partial_tu(t,\cdot)\|_{L^2(M)}^2
\end{equation*}
which is in particular equal to the initial energy $\|\nabla_g u_0\|_{L^2(M)}^2+\|u_1\|_{L^2(M)}^2$. Let $T>0$ and $\omega$ be a measurable subset. The observability inequality reads as follows:
\begin{equation} \label{e:obsRiemwaves}
E(u(0))\leq C\int_0^T\int_\omega |\partial_tu(t,x)|^2dxdt.
\end{equation}
Note that thanks to the conservation of energy, the left hand side of \eqref{e:generobs} has been replaced by the energy of the initial datum.

We set $P=\partial_{tt}^2-\Delta_g$ (which is a second-order pseudo-differential operator), whose principal symbol is
\begin{equation*}
p_2(t,\tau,x,\xi)=-\tau^2+g^*(x,\xi)
\end{equation*}
with $\tau$ the dual variable of $t$ and $g^*$ the principal symbol of $-\Delta_g$. In $T^*(\R\times M)$, the Hamiltonian vector field $\vec{H}_{p_2}$ associated with $p_2$ is given by $\vec{H}_{p_2}f=\{p_2,f\}$. Since $\vec{H}_{p_2}p_2=0$, we get that $p_2$ is constant along the integral curves of $\vec{H}_{p_2}$. Thus, the characteristic set $\mathcal{C}(p_2)=\{p_2=0\}$ is preserved under the flow of $\vec{H}_{p_2}$. Null-bicharacteristics are then defined as the maximal integral curves of $\vec{H}_{p_2}$ which live in $\mathcal{C}(p_2)$.  In other words, the null-bicharacteristics are the maximal solutions of
\begin{equation*}
\left\lbrace \begin{array}{l}
\dot{t}(s)=-2\tau(s)\,,\\
\dot{x}(s)=\nabla_\xi g^*(x(s),\xi(s))\,,\\
\dot{\tau}(s)=0\,,\\
\dot{\xi}(s)=-\nabla_xg^*(x(s),\xi(s))\,,\\
\tau^2(0)=g^*(x(0),\xi(0)).
\end{array}\right.
\end{equation*}
It is well-known that the projection $x(s)$ of a bicharacteristic ray $(x(s),\xi(s))$ traveled at speed $1$ is a geodesic in $M$, i.e., a curve which realizes the minimal distance between any two of its points which are close enough. 

Let us also mention the fact that at the boundary of $M$, the above definition of null-bicharacteristics has to be completed (yielding the so-called Melrose-Sj\"ostrand flow): due to trajectories which ``graze'' along the boundary, one cannot always define the null-bicharacteristics which touch the boundary by reflexion, and we refer the reader to \cite{melrose1978singularities} and \cite{le2017geometric} for more on this subject. In these papers, a notion of ``generalized bicharacteristics'' is defined, which explains how to define bicharacteristics at the boundary. For us, this will only be useful to give a precise statement for Theorem \ref{t:obsRiemwaves}.

\begin{definition}
Let $T>0$ and $\omega\subset M$ be a measurable subset. We say that the Geometric Control Condition holds in time $T$ in $\omega$, and we write $(GCC)_{\omega,T}$, if for any projection $\gamma$ of a bicharacteristic ray traveled at speed $1$, there exists $t\in (0,T)$ such that $\gamma(t)\in\omega$.
\end{definition}
The following result states that the observability of \eqref{e:waveop} is (more or less) equivalent to the geometric condition $(GCC)_{\omega,T}$. It illustrates the finite speed of propagation for waves.
\begin{thm}[\cite{bardos1992sharp}, \cite{burq1997condition}, \cite{humbert2019observability}] \label{t:obsRiemwaves}
Assume that $\omega\neq \emptyset$ is open and that $(GCC)_{\omega,T}$ holds.  Assume also that no generalized bicharacteristic has a contact of infinite order with $(0,T)\times\partial M$. Then \eqref{e:obsRiemwaves} holds, i.e., the wave equation \eqref{e:waveop} is observable in time $T$ on $\omega$. Conversely, if the wave equation \eqref{e:waveop} is observable in time $T$, then $(GCC)_{T,\overline{\omega}}$ holds, where $\overline{\omega}$ denotes the closure of $\omega$.
\end{thm}
Note that the second statement in the last theorem is not the exact converse of the first one, since it involves the closure $\overline{\omega}$ and not simply $\omega$. This is due to the phenomenon of grazing rays: if there exists a ray $\gamma$ which does not enter $\omega$ but which touches the boundary $\partial\omega$, so that the geometric control condition is not satisfied, it can however happen (notably if the flow is ``stable'' close from the ray) that observability holds, see \cite[Section VI.B]{lebeau1992control} for an example.

Considering solutions of \eqref{e:waveop} of the form $e^{it\sqrt\lambda }\varphi$ where $\varphi$ is an eigenfunction of $-\Delta_g$ corresponding to the eigenvalue $\lambda$, the following result follows from Theorem \ref{t:obsRiemwaves}:
\begin{cor} \label{c:eigenfunctions}
Assume that $\omega\neq \emptyset$ is open and that there exists $T>0$ such that $(GCC)_{\omega,T}$ holds. Then, for any eigenfunction $\varphi$ of $-\Delta_g$, there holds
\begin{equation*}
\int_\omega |\varphi(x)|^2 dx \geq C\int_M |\varphi(x)|^2 dx .
\end{equation*}
In particular, $\text{supp}(\varphi)=M$.
\end{cor}
All the observability inequalities stated in this survey paper yield similar lower bounds, but we will not state them thereafter.
\begin{remark}[Gaussian beams] \label{r:gb}
The fact that $(GCC)_{\overline{\omega},T}$ is a necessary condition for observability can be understood as follows. If $(GCC)_{\overline{\omega},T}$ does not hold, then let $\gamma:[0,T]\rightarrow M$ be a geodesic which does not enter $\overline{\omega}$. By compactness, there exists $\varepsilon >0$ such that $\gamma_{|[0,T]}$ does not meet an $\varepsilon$-neighborhood of $\overline{\omega}$. Then, one can construct a sequence of solutions $(u_n)_{n\in\N}$ of the wave equation whose initial energy $E(u_n(0))$ is normalized to $1$, and with energy $E(u(t))$ localized around $\gamma(t)$ at any time $t\in [0,T]$: quantitatively, the energy of $u_n$ outside a tubular neighborhood of $\gamma$ of size $\varepsilon$ tends to $0$ as $n\rightarrow +\infty$. This disproves the observability inequality \eqref{e:obsRiemwaves}. The sequence $(u_n)_{n\in\N}$, if taken as a Gaussian profile centered at a point describing $\gamma$, is called a Gaussian beam.
\end{remark}

\subsubsection{Observability of the Riemannian Schr\"odinger equation} \label{s:obsSchrodRiem}
In case of the Schr\"odinger equation \eqref{e:schrodop}, the observability inequality reads as follows:
\begin{equation} \label{e:obsRiemschrod}
\|u_0\|_{L^2(M)}^2\leq C\int_0^T\int_\omega |u(t,x)|^2dxdt.
\end{equation}
Indeed, as for the wave equation \eqref{e:waveop}, the $L^2$-norm of the solution is preserved along the flow, so that $\|u(T)\|_{L^2}=\|u_0\|_{L^2}$. A sufficient condition for observability is the following:
\begin{thm}[\cite{lebeau1992controle} and Appendix of \cite{dehman2006stabilization}] \label{t:lebschrod}
Assume that $\omega\neq \emptyset$ is open and that $(GCC)_{\omega,T'}$ holds for some $T'>0$. Then \eqref{e:obsRiemschrod} holds, i.e., the Schr\"odinger equation \eqref{e:schrodop} is observable in any time $T>0$ on $\omega$.
\end{thm}
The interplay between $T'$ and $T$ in the above result is due to the fact that the  Schr\"odinger equation ``propagates at infinite speed'' so that no matter how large $T'$ is, observability holds in any time $T>0$ if $(GCC)_{\omega,T'}$ holds. This contrasts with the finite speed of propagation of the wave equation.

The converse of the above theorem, namely to find necessary conditions on $(\omega,T)$ for \eqref{e:obsRiemschrod} to hold, is notoriously a difficult problem. The main results in this direction are for the torus (see \cite{jaffard1990controle}, \cite{burq2012control}, \cite{anantharaman2014semiclassical}), and in Riemannian manifolds with negative curvature (see \cite{dyatlov2019control}), where \eqref{e:obsRiemschrod} holds for any non-empty open subset $\omega$ and any time $T>0$. Indeed, it is expected that if the geodesic flow of the background geometry is unstable, solutions of \eqref{e:schrodop} are more ``delocalized'' than those of \eqref{e:waveop} for example. See also the case of the disk \cite{anantharaman2016wigner}.

\subsubsection{Observability of the Riemannian heat equation}  
Let us end with the heat equation. The observability inequality reads as follows:
\begin{equation} \label{e:obsRiemheat}
\|u(T)\|_{L^2(M)}^2\leq C\int_0^T\int_\omega |u(t,x)|^2dxdt.
\end{equation}
\begin{thm}[\cite{lebeau1995controle}] \label{t:obsRiemheat}
Let $\omega\neq \emptyset$ be open and $T>0$. Then \eqref{e:obsRiemheat} holds, i.e., the heat equation \eqref{e:heatop} is observable in time $T$ on $\omega$.
\end{thm}
Note that no geometric condition on $\omega$ is required in this case. This result illustrates the infinite speed of propagation of the heat equation.

The works presented hereafter address that same problem of observability of linear PDEs, but with focus on subelliptic PDEs, meaning that the Laplace-Beltrami operator is replaced in these PDEs by a sub-Laplacian, which is a subelliptic operator. The next subsection is thus devoted to introduce the main objects of study, namely sub-Laplacians.

\subsection{Sub-Riemannian geometry and sub-Laplacians} After the founding work of Lars H\"ormander, and with the development of sub-Riemannian geometry since the 1980's, subelliptic operators have been considered as a natural generalization of elliptic operators. In particular, sub-Laplacians, which we will define soon, are natural generalizations of the Laplace-Beltrami operator. Therefore, the question of observability/controllability of evolution PDEs driven by sub-Laplacians has been investigated since a decade, with a particular focus on parabolic (or heat-type) equations. In this survey, we mainly focus on wave-type and Schr\"odinger-type subelliptic equations, for which the first results appeared in 2019.

\subsubsection{Sub-Laplacians} \label{s:sL} Instead of defining subelliptic operators in full generality, we prefer here to focus only on sub-Laplacians. The geometry naturally associated to sub-Laplacians is called sub-Riemannian geometry. The books \cite{montgomery2002tour} and \cite{agrachev2019comprehensive} are clear and detailed introductions to this geometry. Readers only interested in the results of Sections \ref{s:sun} and \ref{s:fermanian} could skip this part, and focus on Examples \ref{exGrushin} to \ref{exHtype} which are sufficient for these sections. 

\smallskip

Let $n\in\mathbb{N}^*$ and let $M$ be a smooth connected compact manifold of dimension $n$ with a non-empty boundary $\partial M$. Let $\mu$ be a smooth volume on $M$.

\smallskip

We consider $m\geq 1$ smooth vector fields $X_1,\ldots,X_m$ on $M$ which are not necessarily independent, and we assume that the following Lie bracket generating (or H\"ormander) condition holds (see \cite{hormander1967hypoelliptic}):

\begin{center}
The vector fields $X_1,\ldots,X_m$ and their iterated brackets $[X_i,X_j], [X_i,[X_j,X_k]]$, etc. span the tangent space $T_qM$ at every point $q\in M$.
 \end{center}

\smallskip
We consider the sub-Laplacian $\Delta$ defined by
\begin{equation} \label{d:subLapl}
\Delta=-\sum_{i=1}^m X_i^*X_i=\sum_{i=1}^m X_i^2+\div_\mu(X_i)X_i
\end{equation}
where the star designates the transpose in $L^2(M,\mu)$ and the divergence with respect to $\mu$ is defined by $L_X\mu=(\div_\mu X)\mu$, where $L_X$ stands for the Lie derivative. Then $\Delta$ is hypoelliptic (see \cite[Theorem 1.1]{hormander1967hypoelliptic}). 

We set
$$
\mathcal{D}={\rm Span}(X_1,\ldots,X_m)\subset TM
$$
which is called the \emph{distribution} associated to the vector fields $X_1,\ldots,X_m$. For $x\in M$, we denote by $\mathcal{D}_x$ the distribution $\mathcal{D}$ taken at point $x$. Note that $\mathcal{D}$ does not necessarily have constant rank. When  $\mathcal{D}=TM$, the operator $\Delta$ is elliptic.

\smallskip

We also introduce the \emph{metric} $g$ on $\mathcal{D}$ defined at any $x\in M$ by
$$
g_x(v,v)=\inf\left\{\sum_{i=1}^m u_i^2 \ | \ v=\sum_{i=1}^m u_iX_i(x) \right\}.
$$
This is a Riemannian metric on $\mathcal{D}$. We call $(M,\mathcal{D},g)$ a \emph{sub-Riemannian structure}.

\smallskip

 In the general case where $\mathcal{D}\subsetneq TM$, the set $TM\setminus \mathcal{D}$ can be understood as the directions where the metric $g$ takes the value $+\infty$. A well-known theorem, due to Chow and Rashevskii, asserts that any two points can be joined by a path, i.e., a continuous function $\gamma:[0,1]\rightarrow M$ with derivative $\dot\gamma(t)$ contained in $\mathcal{D}_{\gamma(t)}$ for almost any $t\in[0,1]$. In other words, the sub-Riemannian distance
\begin{equation*}
d_g(x_0,x_1)=\inf\left\{\int_0^1\sqrt{g_{\gamma(t)}(\dot{\gamma}(t),\dot{\gamma}(t))}dt \ \big| \ \gamma(0)=x_0, \gamma(1)=x_1, \ \dot{\gamma}(t)\in\mathcal{D}_{\gamma(t)} \text{ a.s. for } t\in [0,1]\right\}
\end{equation*}
is finite for any $x_0,x_1\in M$.

When moving in a sub-Riemannian structure, $\mathcal{D}$ should be understood as the ``set of allowed directions for the motion'', and, although it is not possible to move directly in directions of $TM\setminus \mathcal{D}$, Chow-Rashevskii's theorem asserts that any two points can be joined by a path. This is due to ``indirect motions'', that is, paths which describe spirals turning around a fixed forbidden direction of $TM\setminus \mathcal{D}$ and thus advancing  in this direction (although indirectly).

\begin{definition} \label{d:step}
The step $k$ is the least integer $k\in\N$ such that the Lie brackets of the vector fields $X_1,\ldots,X_m$ of length $\leq k$  (i.e., $X_i$, $[X_i, X_j]$, $[X_i,[X_j,X_\ell]]$, up to length $k$)  span the whole tangent space $TM$.
\end{definition}
\begin{remark} \label{r:step}
More generally, the step $k_x$ can be defined at any point $x\in M$, just by considering the Lie brackets of the vector fields $X_1,\ldots,X_m$ at point $x$.
\end{remark}


\subsubsection{Examples} \label{s:examples} We now give a few examples of sub-Laplacians which we shall study in the sequel. 
\begin{example} \label{exGrushin}
On $M=\R_x\times\R_y$, we set $\Delta_G=\partial_x^2+x^2\partial_y^2$. This sub-Laplacian is the so-called Baouendi-Grushin operator, sometimes unproperly called simply Grushin operator (see \cite[Section 11]{garofalo2017fractional}). In this case, $X_1=\partial_x$, $X_2=x\partial_y$ and ${\rm Span}(X_1,X_2,[X_1,X_2])=TM$. In particular, $\mathcal{D}=TM$ outside the line $\{x=0\}$. Also, $\mu$ is the Lebesgue measure. The step is $2$ on the line $\{x=0\}$ and $1$ outside this line. Since the sub-Riemannian structure is ``Riemannian'' outside this line, the Baouendi-Grushin operator is sometimes called ``almost-Riemannian''.
\end{example}
\begin{example} \label{exgamma}
More generally, for $\gamma\geq 0$ (not necessarily an integer), one can consider $\Delta_\gamma=\partial_x^2+|x|^{2\gamma}\partial_y^2$ on the same manifold $M=(-1,1)_x\times\T_y$. For $\gamma\in \N$, the step is $k=\gamma+1$. Note that for $\gamma\notin \N$, the H\"ormander condition is not necessarily satisfied, but we include this class of examples in our study since our computations allow to handle them.
\end{example}
\begin{example} \label{exHeis}
Given $d\in\N^*$, one can also define a sub-Laplacian arising from the Heisenberg group $\mathbf{H}^d$ of dimension $2d+1$. Recall that the Heisenberg group~$\mathbf{H}^d$ is $\R^{2d+1}$ endowed with the group law $(x,y,z) \cdot (x',y',z'):=(x+x', y+y',z+z'+\frac12 \sum_{j=1}^d(x_jy_j'-x_j'y_j))$,
where $x,y,x',y'\in\R^d$ and $z,z'\in\R$. Taking the left-quotient of $\mathbf{H}^d$ by the discrete subgroup $\widetilde{\Gamma}=(\sqrt{2\pi}\mathbb{Z})^{2d}\times \pi\mathbb{Z}$, we obtain a compact manifold $M=\widetilde{\Gamma}\backslash \mathbf{H}^d$. Let
\begin{equation*}
X_j=\partial_{x_j}-\frac{y_j}{2}\partial_{z}, \ \  \ Y_j=\partial_{y_j}+\frac{x_j}{2}\partial_{z}, \qquad \text{for } j=1,\ldots,m,
\end{equation*}
which are left-invariant and can be thus considered as vector fields on the quotient manifold $M$. Finally, we define the sub-Lapacian
\begin{equation*}
\Delta_{\widetilde{\Gamma}\backslash \textbf{H}^d}=\sum_{j=1}^d X_j^2+Y_j^2.
\end{equation*}
Since $[X_j,Y_j]=\partial_z$ for any $j$, the step is $2$.
\end{example}
\begin{example} \label{exHtype}
Heisenberg-type groups are generalizations of Heisenberg groups.  They were first introduced in \cite{kaplan1980fundamental}, where the fundamental solution of the associated sub-Laplacians, which is particularly simple, was computed. These groups give rise to step $2$ sub-Riemannian structures, which center can be of dimension $p>1$ (whereas the center of the Heisenberg group $\mathbf{H}^d$ is of dimension $1$). For a precise definition, see \cite{kaplan1980fundamental}.
\end{example}
\begin{example} \label{excontact}
Contact sub-Laplacians are associated to sub-Riemannian structures of ``contact type''. We assume that the vector fields $X_1,\ldots,X_m$ span a distribution $\mathcal{D}$ which is a contact distribution over $M$, i.e., $M$ has odd dimension $n=2m+1$ and there exists a $1$-form $\alpha$ on $M$ with $\mathcal{D}={\rm Ker}(\alpha)$ and $\alpha\wedge (d\alpha)^m\neq 0$ at any point of $M$. Then, for any smooth volume $\mu$, the sub-Laplacian $\Delta$ is called a contact sub-Laplacian. A typical example is given by the Heisenberg sub-Laplacian $\Delta_{\widetilde{\Gamma}\backslash \textbf{H}^d}$ defined above. 
\end{example}
\begin{example} \label{exanormales}
Magnetic Laplacians are also sub-Laplacians, and we focus here on a simple family of examples. In $\R^3$ with coordinates $x,y,z$, we consider the two vector fields $X_1=\partial_x-A_x(x,y)\partial_z$ and $X_2=\partial_y-A_y(x,y)\partial_z$ where $A_x$, $A_y$ are functions which do not depend on $z$. The magnetic Laplacian is then $\Delta=X_1^2+X_2^2$. The $1$-form $A=A_xdx+A_ydy$ is called the connection form, and the $2$-form $B=dA$ is called the magnetic field. The modulus $|b|$ of the function $b$ defined by the relation $B=b\;dx\wedge dy$ is called the intensity of the magnetic field. Taking a quotient or assuming that $|b|$ is bounded away from $0$ at infinity, it is possible to assume that the sub-Laplacian $\Delta$ has a compact resolvent.

Magnetic Laplacians were used by Montgomery to prove the existence of abnormal minimizers in some sub-Riemannian geometries (see \cite{montgomery1994abnormal}): abnormal minimizers are local minimizers of the sub-Riemannian distance which are not projections of bicharacteristics, and they show up for example as zero curves of the intensity $b$. Subsequently, Montgomery showed that the spectral asymptotics of $\Delta$ are very different depending on the fact that $b$ vanishes or not (see  \cite{montgomery1995hearing}).
\end{example}

\begin{remark}
In all the previous examples, as well as in the general definition of sub-Laplacians given above, the volume $\mu$ is assumed to be smooth. However, one can also define sub-Laplacians for non-smooth volumes $\mu$. This is natural when the sub-Riemannian distribution is singular, for example $\mathcal{D}={\rm Span}(\partial_x,x\partial_y)$ in $\R^2$ (see Example \ref{exGrushin}), since the Popp measure, which is an intrinsic measure defined on sub-Riemannian manifolds (see \cite[Section 10.6]{montgomery2002tour}), ``blows up''. The associated sub-Laplacian is then unitarily equivalent to a sub-Laplacian with smooth volume plus a singular potential. The essential self-adjointness of some of these sub-Laplacians has been studied for example in \cite{boscain2013laplace}. When they are essentially self-adjoint, the controllability/observability of the associated evolution equations can also be studied: this is an open question which we do not address here.
\end{remark}

\subsubsection{Hypoellipticity and subellipticity} \label{s:hyposub} Two notions are often used to qualify the smoothing properties of sub-Laplacians: hypoellipticity and subellipticity. Here, we briefly recall their definitions and explain why they are not exactly equivalent.

\begin{definition} 
A (pseudo-)differential operator $A$ with $C^\infty$ coefficients in $M$ is \emph{hypoelliptic} in $M$ if for all $u\in\mathcal{D}'(M)$ and $x\in M$, if $Au\in C^\infty$ near $x$, then $u\in C^\infty$ near $x$.
\end{definition}
Hypoellipticity appeared naturally in the work of Kolmogorov \cite{kolmogoroff1934zufallige} on the motion of colliding particles when he wrote down the equation
\begin{equation*}
\partial_t u -\mathcal{L}u=f \qquad \text{where} \qquad \mathcal{L}=x\partial_y+\partial_x^2.
\end{equation*}
Indeed, the operator $\mathcal{L}$ is hypoelliptic.

\begin{definition} 
A formally selfadjoint (pseudo-)differential operator $A:C^\infty(M)\rightarrow C^\infty(M)$ of order $2$ is said to be \emph{subelliptic} if there exist $s,C>0$ such that 
\begin{equation}  \label{e:subest}
\|u\|_{H^s(M)}^2\leq C((Au,u)_{L^2(M)}+\|u\|_{L^2(M)}^2)
\end{equation}
 for any $u\in C^\infty(M)$.
 \end{definition}

Under the Lie bracket condition (H\"ormander condition), H\"ormander was able to prove that any sub-Laplacian $\Delta$ is hypoelliptic (see \cite{hormander1967hypoelliptic} and \cite[Chapter 2]{helffer2005hypoelliptic}).  His proof relies on the fact that $\Delta$ is subelliptic; indeed, the optimal $s$ in \eqref{e:subest} is $1/k$, where $k$ is the step of the associated sub-Riemannian structure, as proved by Rotschild and Stein \cite[Theorem 17 and estimate (17.20)]{rothschild1976hypoelliptic}.

Conversely, note that an hypoelliptic ``sum of squares'' (i.e., an operator of the form \eqref{d:subLapl} which is hypoelliptic) does not necessarily satisfy the Lie bracket assumption : given a smooth function $a:\R\rightarrow\R$ vanishing at infinite order at $0$ but with $a(s)>0$ for $s\neq0$,  the sub-Laplacian $\Delta=\partial_{x_1}^2+a(x_1)^2\partial_{x_2}^2$ on $\R_{x_1x_2}^2$ is hypoelliptic although the Lie bracket condition fails (see \cite{fedii1971criterion} and \cite{morimoto1978hypoellipticity}).

Let us finally mention that some operators $A$ satisfy the property that if $Au$ is real-analytic, then $u$ is real-analytic: they are called {\it analytic hypoelliptic}. The so-called Tr\`eves conjecture describes a possible link between analytic hypoellipticity of an operator and the absence of abnormal geodesics (see \cite{treves1999symplectic} for the conjecture and \cite{albano2018analytic} for more recent results).

We end this section with a remark concerning the compactness of the manifolds $M$ considered in this survey.
\begin{rem}
Because of the physical nature of the problems studied in control/observability theory, most equations are set in compact manifolds, and this survey is no exception to the rule. Even in Example \ref{exHeis}, the sub-Laplacian is defined on a compact quotient of the Heisenberg group. Together with the hypoellipticity, the compactness of the underlying manifold implies that all sub-Laplacians have a compact resolvent, and thus a discrete spectrum, which is of importance for deriving properties of eigenfunctions from observability results.
\end{rem}

\subsection{Observability of subelliptic PDEs: known results}
This section is devoted to stating results which were previously known in the literature about controllability/observability of subelliptic PDEs. All PDEs we consider are well-posed in natural energy spaces which we do not systematically recall. 

\subsubsection{Subelliptic heat equations}  \label{s:subheatequations} 
Let us start with the result proved in \cite{beauchard2014null}, which concerns the heat equation \eqref{e:heatop} where $\mathcal{L}=\Delta_\gamma$ is the Baouendi-Grushin-type operator introduced in Example \ref{exgamma} for some $\gamma>0$. The open subset of observation $\omega\subset (-1,1)\times \T$ they consider is a vertical strip of the form $(a,b)\times \T$ where $0<a<b<1$. The observability inequality is \eqref{e:obsRiemheat}, with the modification that $u$ runs over the set of solutions of \eqref{e:heatop} with $\mathcal{L}=\Delta_\gamma$. The authors prove the following result, to be compared with Theorem \ref{t:obsRiemheat}:
\begin{thm}[\cite{beauchard2014null}] \label{t:bcg14}
Let $\gamma>0$ and $\omega$ be as above. Then
\begin{itemize}
\item If $\gamma\in (0,1)$, then for any $T>0$, \eqref{e:obsRiemheat} holds;
\item If $\gamma=1$, i.e., $\Delta_\gamma=\Delta_G$, then there exists $T_0>0$ such that \eqref{e:obsRiemheat} holds if $T>T_0$ and does not hold if $T<T_0$;
\item If $\gamma>1$, then, for any $T>0$, \eqref{e:obsRiemheat} fails.
\end{itemize}
\end{thm}
The proof is done by establishing an infinite number of Carleman inequalities for operators $-\partial_x^2+n^2x^2$ for $n\in\Z$, with bounds uniform in $n$. It was proved  in \cite{beauchard2015grushin} that in case $\gamma=1$, the minimal time $T_0$ is equal to $a^2/2$. The fact that $T_0\geq a^2/2$ can be seen by using explicit eigenfunctions of $\Delta_\gamma$.

Koenig studied the observability of \eqref{e:heatop} with $\mathcal{L}=\Delta_G$, but for another geometry of the observation set $\omega$: this time, it is a horizontal band of the form $(-1,1)\times I$ where $I$ is a proper open subset of $\T$. 
\begin{thm}[\cite{koenig2017non}] \label{t:koenig}
Let $\omega=(-1,1)\times I$ where $I$ is a proper open subset of $\T$. Then \eqref{e:obsRiemheat} fails for any $T>0$.
\end{thm} 
The proof of this result relies on the non-observability of the 1D half-heat equation $\partial_tu+|\partial_x|u=0$ and on techniques coming from complex analysis where the complex variable is $z=e^{-t+iy}$. 

Although the observability properties of the heat equation driven by general hypoelliptic operators are still mysterious, we list here a few works addressing this question. The recent works \cite{lissy2020non}, \cite{beauchard2020minimal} and \cite{duprez2020control} continue and generalize the analysis of \cite{beauchard2014null} and \cite{koenig2017non} on the control of the Baouendi-Grushin heat equation. Besides, \cite{beauchard2017heat} establishes the existence of a minimal time of observability, as in the second point of Theorem \ref{t:bcg14}, for the heat equation driven by the Heisenberg sub-Laplacian of Example \ref{exHeis}. Let us finally mention the papers \cite{darde2020critical} and \cite{beauchard2018null} which also deal with controllability issues for hypoelliptic parabolic equations. 

The above theorems show that some subelliptic heat equations driven by simple sub-Laplacians require a larger time to be observable than the usual Riemannian heat equation, and observability may even fail in any time $T>0$. As we will see, this is a very general phenomenon for subelliptic evolution PDEs, at least for subelliptic wave equations and (some) Schr\"odinger-type equations. Our results, however, do not shed any new light on subelliptic heat equations, which remain mysterious due to the lack of ``general arguments'' which would not rely on geometric and analytic features specific to very particular sub-Laplacians.

\subsubsection{Approximate observability of subelliptic PDEs}
Recently, Laurent and L\'eautaud have studied the observability of subelliptic PDEs but with focus on a different notion of observability, called approximate observability. As before, all their results can also be stated in terms of a dual notion, called approximate controllability, however we will not even mention these dual statements in order to keep the presentation as simple as possible. Their results are quantitative, in the sense that they give explicit bounds on the observability constants involved in their results. The next paragraphs are devoted to a brief description of their results (see \cite{laurent2020tunneling}).

Let us consider a sub-Laplacian $\Delta$ as in \eqref{d:subLapl}, with associated sub-Riemannian structure $(M,\mathcal{D},g)$. We assume that the manifold $M$ (assumed to have no boundary, $\partial M=\emptyset$), the smooth volume $\mu$ and the vector fields $X_i$ are all {\it real-analytic}. For $s\in\R$, the operator $(1-\Delta)^{\frac{\ell}{2}}$ is defined thanks to functional calculus, and we consider the (adapted) Sobolev spaces
\begin{equation*}
\mathcal{H}^\ell(M)=\{u\in\mathcal{D}'(M), (1-\Delta)^{\frac{\ell}{2}}u\in L^2(M)\}
\end{equation*}
with the associated norm $\|u\|_{\mathcal{H}^\ell(M)}=\|(1-\Delta)^{\ell}u\|_{L^2(M)}$.

\begin{thm}[\cite{laurent2020tunneling}] \label{t:LL20}
Let $\omega$ be a non-empty open subset of $M$ and let $T>\sup_{x\in M}d_g(x,\omega)$. We denote by $k$ the step. Then there exist $c,C>0$ such that
\begin{equation} \label{e:weakobsLL20}
\|(u_0,u_1)\|_{\mathcal{H}^{1}\times L^2}\leq Ce^{c\Lambda^k}\|u\|_{L^2((-T,T)\times \omega)}, \qquad \text{with  } \Lambda=\frac{\|(u_0,u_1)\|_{\mathcal{H}^{1}\times L^2}}{\|(u_0,u_1)\|_{L^2\times \mathcal{H}^{-1}}}
\end{equation}
for any solution $u$ of \eqref{e:waveop} on $(-T,T)$ such that $(u,\partial_tu)_{|t=0}=(u_0,u_1)\in\mathcal{H}^1(M)\times L^2(M)$. 
\end{thm}
The above result in particular implies unique continuation (and quantifies it): if $u= 0$ in $(-T,T)\times \omega$, then $u\equiv 0$. However, the exact observability inequality which we shall study (see \eqref{e:obsRiemwaves}) is a stronger requirement than \eqref{e:weakobsLL20}, in particular because of the presence of the ``typical frequency of the datum'' $\Lambda$ in the right-hand side of \eqref{e:weakobsLL20}. The techniques used for proving Theorem \ref{t:LL20} are totally different from those we present in the sequel.

\subsubsection{Observability of Baouendi-Grushin Schr\"odinger equation} The recent work \cite{burq2019time} is the first one dealing with exact observability of a subelliptic  Schr\"odinger equation, namely in the context of Example \ref{exGrushin} with observation set given by a horizontal band as in Theorem \ref{t:koenig}. The observability inequality is given by \eqref{e:obsRiemschrod}, except that $u$ runs over the solutions of the Schr\"odinger equation driven by the sub-Laplacian $\Delta_G$.
\begin{thm}[\cite{burq2019time}] \label{t:burqsun}
Let $M=(-1,1)\times\T$ and $\Delta_G=\partial_x^2+x^2\partial_y^2$. Let $\omega=(-1,1)\times I$ where $I\subsetneq \T$ is open. Let $T_0=\mathscr{L}(\omega)$ be the length of the maximal sub-interval contained in $\T\setminus I$. Then, the observability property \eqref{e:obsRiemschrod} holds if and only if $T> T_0$.
\end{thm}
Again, this result shows the existence of a minimal time of control which contrasts with the ``infinite speed of propagation'' illustrated by Theorem \ref{t:lebschrod}. Its proof relies on fine semi-classical analysis, somehow linked to that explained in Section \ref{s:sun}.

\subsubsection{Non-linear subelliptic PDEs} Although this survey is devoted only to \emph{linear} subelliptic PDEs, let us say a word about \emph{non-linear} subelliptic PDEs. To study the cubic Grushin-Schr\"odinger equation $i\partial_tu-(\partial_x^2+x^2\partial_y^2)u=|u|^2u$, Patrick G\'erard and Sandrine Grellier introduced a toy model, the cubic Szeg\"o equation, which models the interactions between the nonlinearity and the lack of dispersivity of the linear equation (already visible in the above Theorem \ref{t:burqsun}). In \cite{gerard2010cubic}, they put this equation into a Hamiltonian framework and classify the traveling waves for this equation, which are related to the traveling waves resulting from \eqref{e:decompoL2grushin}.

\subsection{Main results}
Let us now present the main results contained in the papers \cite{letrouit2020subelliptic}, \cite{letrouit2020observability} and \cite{kammerer2020observability}. All of them illustrate the slowdown of propagation of evolution PDEs in directions transverse to the distribution: in a nutshell, observability will require a much longer time to hold for subelliptic PDEs than for elliptic ones, and this time will be even larger when the step $k$ is larger. All our results are summarized in Figure \ref{t:summary} at the end of this section.

\subsubsection{First main result} We start with a general result on subelliptic wave equations. Let $\Delta=-\sum_{i=1}^m X_i^*X_i$ be a sub-Laplacian, where the adjoint denoted by star is taken with respect to a volume $\mu$ on $M$, which is assumed to have a boundary $\partial M\neq \emptyset$.\footnote{This assumption is not necessary, since Theorem \ref{t:main1} also works for manifolds without boundary, but this would require to introduce a slightly different notion of observability.} Consider the wave equation
\begin{equation} \label{e:system}
\left\lbrace \begin{array}{l}
\partial_{tt}^2u-\Delta u=0\, \text{ \ in $(0,T)\times M$}\\
u=0 \text{ \ on } (0,T)\times \partial M, \\
(u_{|t=0},\partial_tu_{|t=0})=(u_0,u_1)\,
\end{array}\right.
\end{equation}
where $T>0$, and the initial data $(u_0,u_1)$ are in an appropriate energy space.
 The natural energy of a solution $u$ of the sub-Riemannian wave equation \eqref{e:system} is
 \begin{equation*}
 E(u(t,\cdot))= \frac12\int_M \left( |\partial_tu(t,x)|^2 +\sum_{j=1}^m (X_ju(t,x))^2 \right) d\mu(x).
 \end{equation*}
 Observability holds in time $T_0$ on $\omega$ if there exists $C>0$ such that for any solution $u$ of \eqref{e:system}, 
 \begin{equation} \label{e:strongobs}
E(u(0))\leq C \int_0^{T_0} \int_\omega|\partial_tu(t,x)|^2d\mu(x)dt.
 \end{equation}
\begin{thm}[\cite{letrouit2020subelliptic}] \label{t:main1}
Let $T_0>0$ and let $\omega\subset M$ be a measurable subset. We assume that there exist $1\leq i,j\leq m$ and $x$ in the interior of $M\backslash\omega$ such that $[X_i,X_j](x)\notin {\rm Span}(X_1(x),\ldots,X_m(x))$. Then the subelliptic wave equation \eqref{e:system} is not exactly observable on $\omega$ in time $T_0$.
\end{thm}
Theorem \ref{t:main1} can be reformulated as follows: subelliptic wave equations are never observable. The condition that there exist $1\leq i,j\leq m$ and $x$ in the interior of $M\backslash\omega$ such that $[X_i,X_j](x)\notin {\rm Span}(X_1(x),\ldots,X_m(x))$ means that $\Delta$ is not elliptic at $x$; this  assumption is absolutely necessary since otherwise, locally, \eqref{e:system} would be the usual elliptic wave equation, and its observability properties would depend on the GCC, as stated in Theorem \ref{t:obsRiemwaves}. 

The key ingredient in the proof of Theorem \ref{t:main1} is that the GCC fails for any time $T_0>0$: in other words, there exist geodesics which spend a time greater than $T_0$ outside $\omega$. Then, the Gaussian beam construction described in Remark \ref{r:gb} allows to contradict the observability inequality \eqref{e:strongobs}.

\subsubsection{Second main result} Our second main result, obtained in collaboration with Chenmin Sun, sheds a different light on Theorem \ref{t:main1}. For this second statement, we consider the generalized Baouendi-Grushin operator of example \ref{exgamma}, i.e, $\Delta_\gamma=\partial_x^2+|x|^{2\gamma}\partial_y^2$ on $M=(-1,1)_x\times\T_y$, and we assume that $\gamma\geq 1$ (not necessarily an integer). We also consider the Schr\"odinger-type equation with Dirichlet boundary conditions
\begin{equation} \label{e:schrodfrac}
\left\lbrace \begin{array}{l}
i\partial_tu-(-\Delta_\gamma)^{s}u=0 \\
u_{|t=0}=u_0 \in L^2(M) \\
u_{|x=\pm 1}=0
\end{array}\right.
\end{equation}
where $s\in\N$ is a fixed integer. 
Given an open subset $\omega\subset M$, we say that \eqref{e:schrodfrac} is observable in time $T_0>0$ in $\omega$ if there exists $C>0$ such that for any $u_0\in L^2(M)$, 
\begin{equation} \label{e:obsschrod}
\|u_0\|_{L^2(M)}^2\leq C\int_0^{T_0} \|e^{-it(-\Delta_\gamma)^s}u_0\|_{L^2(\omega)}^2 dt.
\end{equation}
Our second main result, obtained in collaboration with Chenmin Sun, roughly says that observability holds if and only if the subellipticity (measured by the step $\gamma+1$ in case $\gamma\in\mathbb{N}$), is not too strong compared to the strength of propagation $s$:
\begin{thm}[\cite{letrouit2020observability}] \label{t:main2}
Assume that $\gamma\geq 1$. Let $I\subsetneq \T_y$ be a strict open subset, and let $\omega=(-1,1)_x\times I$. Then, for $s\in\N$, we have:
\begin{enumerate}
\item If $\frac12(\gamma+1)<s$, \eqref{e:schrodfrac} is observable in $\omega$ for any $T_0>0$;
\item If  $\frac12(\gamma+1)=s$, there exists $T_{\inf}>0$ such that \eqref{e:schrodfrac} is observable in $\omega$ for $T_0$ if and only if $T_0> T_{\inf}$;
\item If $\frac12(\gamma+1)>s$, for any $T_0>0$, \eqref{e:schrodfrac} is not observable in $\omega$.
\end{enumerate}
\end{thm}
The case $s=1/2$ corresponds to wave equations. Strictly speaking, it is not covered by Theorem \ref{t:main2} since $s$ is assumed to belong to $\N$ in this theorem, but we see that for any positive $\gamma$ it is roughly related to Point (3), and we thus recover the intuition given by Theorem \ref{t:main1} that subelliptic wave equations should not be observable. The case $\gamma=s=1$ allows to recover Theorem \ref{t:burqsun}, except that we do not find with our method the critical time $T_{\inf}$. Let us also notice that if $\gamma\in\N$, since $\gamma+1$ is the step of the sub-Laplacian $\Delta_\gamma$, the number $\frac12(\gamma+1)$ appearing in Theorem \ref{t:main2} coincides with the exponent known as the gain of Sobolev derivatives in subelliptic estimates (see Section \ref{s:hyposub}). 

\subsubsection{Third main result} \label{s:thirdmain} Finally, our third main result, obtained in collaboration with Clotilde Fermanian Kammerer, illustrates how tools coming from noncommutative harmonic analysis can be used to analyze sub-Laplacians and the associated evolution equations. Our main message is that \emph{a pseudodifferential calculus ``adapted to the sub-Laplacian'' can be used to prove controllability and observability results for subelliptic PDEs} (instead of the usual pseudodifferential calculus used for example to prove Theorem \ref{t:main2}). As we will see, in the present context, once defined this natural pseudodifferential calculus and the associated semi-classical measures (which relies essentially on functional analysis arguments), observability results follow quite directly. 

To relate this last result to the previous ones, let us say that it is roughly linked to the critical case $s=\gamma=1$ of Point 2 of Theorem \ref{t:main2}, i.e., to the case where subelliptic effects are exactly balanced by the strength of propagation of the equation. Indeed, we consider the usual Schr\"odinger equation ($s=1$) in some particular non-commutative Lie groups, called H-type groups (see Example \ref{exHtype}), which have step $2$ (corresponding to $\gamma=1$ for Baouendi-Grushin operators). As in Point 2 of Theorem \ref{t:main2}, we establish that under some geometric conditions on the set of observation $\omega$, observability holds if and only if time is sufficiently large. The main difference with Theorem \ref{t:main2} relies in the tools used for the proof, which could lead to different generalizations. For example, the tools employed in this section allow to handle the case with analytic potential, see \eqref{e:Schrod} below. Also, with these tools, we could imagine to prove observability results for higher-step nilpotent Lie groups, but it requires to know explicit formulas for their representations, since they determine the propagation properties of the semi-classical measures we construct (see Proposition \ref{p:measure0}).

To keep the presentation as simple as possible, we will present our last result only for the Heisenberg groups $\mathbf{H}^d$ of Example \ref{exHeis}, and not for general H-type groups (which are handled in \cite{kammerer2020observability}). By doing so, we avoid defining general H-type groups, while keeping the main message of this work, namely the use of noncommutative harmonic analysis for proving observability inequalities.

Using the notations of Example \ref{exHeis}, we consider $M=\widetilde{\Gamma}\backslash \textbf{H}^d$ together with the equation
  \begin{equation}  \label{e:Schrod}
i\partial_t u+\frac 12 \Delta_M u  +\mathbb{V} u =0
\end{equation}
on $M$, where $\mathbb{V}$ is an analytic function defined on~$M$. The factor $\frac12$ in front of $\Delta_M$ plays no role, we put it here just to keep the same conventions as in \cite{kammerer2020observability}.

The Schr\"odinger equation \eqref{e:Schrod} is observable in time $T$ on the measurable set $U$ if there exists a constant $C_{T,U}>0$ such that
\begin{equation}\label{obs}
\forall u_0\in L^2(M),\;\;\|u_0\|^2_{L^2(M)} \leq C_{T,U} \int_0^T \left\|  {\rm e}^{it(\frac12 \Delta_M+\mathbb{V})} u_0\right\|^2_{L^2(U)} dt.
\end{equation}
Recall that Theorem \ref{t:lebschrod} asserts that, in the Riemannian setting and without potential, the observability of the Schr\"odinger equation is implied by the Geometric Control Condition (GCC), which says that any trajectory of the geodesic flow enters $U$ within time $T$. Using normal geodesics, one can also define a sub-Riemannian geodesic flow (see Section \ref{s:exactgb}) but in some directions of the phase space (called degenerate directions in the sequel), its trajectories are stationary.  For them, we thus need to replace GCC by another condition.  In the case of the Heisenberg group $\mathbf{H}^d$, there is only one such direction, thought as ``vertical'' since it is related to the $\partial_z$ vector field. 

The Heisenberg group $\mathbf{H}^d$ comes with a Lie algebra~${\mathfrak g}$. Via the exponential map
 $$
{\rm Exp} :  {\mathfrak g} \rightarrow \mathbf{H}^d
 $$
 which is a diffeomorphism from ${\mathfrak g}$ to $\mathbf{H}^d$,
 one identifies $\mathbf{H}^d$ and ${\mathfrak g}$ as a set and a manifold. Moreover, ${\mathfrak g}$  is  equipped with a vector space decomposition
$$
  \displaystyle
  {\mathfrak g}=  \mathfrak v \oplus \mathfrak z \, ,
  $$
  such that $[{\mathfrak v},{\mathfrak v}]= {\mathfrak z}\not=\{0\}$ and ${\mathfrak z}$ (of dimension $1$) is the center of ${\mathfrak g}$. We define a scalar product on $\mathfrak z$ by saying that $\partial_z$ has norm $1$, which allows to identify $\mathfrak z$ to its dual $\mathfrak z^*$. We also fix an orthonormal basis $V=(V_1,\ldots, V_{2d})$ of $\mathfrak v$.

We consider the ``vertical'' flow map (also called ``Reeb'', in honor of Georges Reeb) on~$M\times \mathfrak z^*$:
$$\Phi^s_0: (x,\lambda)\mapsto ({\rm Exp} (sd\mathcal Z^{(\lambda)}/2)x,\lambda),\qquad s\in\R$$
 where, for $\lambda\in\mathfrak z$,  $\mathcal Z^{(\lambda)}$ is the element of~$\mathfrak z$ defined by $\lambda(\mathcal Z^{(\lambda)})=|\lambda|$ (or equivalently, $\mathcal Z^{(\lambda)}= \lambda/|\lambda|$ after identification of $\mathfrak z$ and $\mathfrak z^*$). We introduce the following H-type geometric control condition.
{\bf (H-GCC)} The measurable set $U$ satisfies  \textbf{H-type GCC}  in time $T$ if
$$\forall  (x,\lambda)\in M\times( \mathfrak z^*\setminus\{0\}),\;\; \exists s\in (0,T),\;\;\Phi^{s}_0 ((x,\lambda))\in U\times \mathfrak z^*.$$
The flow $\Phi_0^s$ thus replaces the geodesic flow in the degenerate direction.
\begin{definition}
We denote by $T_{\rm GCC}(U)$ the infimum of all $T>0$ such that H-type GCC holds in time~$T$ (and we set $T_{\rm GCC}(U)=+\infty$ if H-type GCC does not hold in any time).
\end{definition}
We also consider the additional  assumption:
\begin{itemize}
\item[{\bf (A)}]  The lift in $\mathbf{H}^d$ of any geodesic of $\mathbf{T}^{2d}$ enters $\omega$ in finite time. \footnote{This condition can be more precisely stated as follows. For any $(x,\omega)\in M\times \mathfrak v^*$ such that $|\omega|=1$, there exists $s\in\R$ such that ${\rm Exp} (s\omega\cdot V) x \in U$. Here, ${\omega\cdot V=\sum_{j=1}^{2d} \omega_jV_j}$ where $\omega_j$ denote the coordinates of $\omega$ in the dual basis of $V$ and it is assumed that ${\sum_{j=1}^{ 2d} \omega^2_j=1}$.}
\end{itemize}

\begin{thm}[\cite{kammerer2020observability}]\label{t:main3} Let $U\subset M$ be open and denote by $\overline{U}$ its closure.
\begin{enumerate}
\item Assume that $U$ satisfies {\bf (A)} and that $T> T_{\rm GCC}(U)$, then the observability inequality \eqref{obs} holds, i.e. the Schr\"odinger equation~\eqref{e:Schrod} is observable  in time $T$ on~$U$.
\item Assume $T\leq T_{\rm GCC}(\overline{U})$, then the observability inequality \eqref{obs} fails.
\end{enumerate}
\end{thm}
This statement looks like Theorem \ref{t:obsRiemwaves} which holds for elliptic waves. In some sense, ``the Schr\"odinger equation in Heisenberg groups looks like an elliptic wave equation'', a phenomenon which was already pointed out by authors studying Strichartz estimates, see \cite{bahouri2000espaces} and \cite{bahouri2016dispersive} for example.

Let us also say that, as already mentioned, Theorem \ref{t:main3} holds more generally in quotients of H-type groups.

To conclude, let us draw a table summing up most of the results presented in this introduction:
 \begin{center}
 \begin{figure}[h]
\begin{tabular}{|c|c|c|c|c|}
  \hline
  & Elliptic & Step 2 & Step $4$ & Step $>4$\\
  \hline
  Waves and half-waves ($s=1/2$) & $T_{\rm inf}$ \footnotesize(under GCC) & $\infty$&$\infty$&$\infty$\\
  \hline
  Schr\"odinger ($s=1$) & 0 \footnotesize(under GCC) & \color{blue}{$T_{\rm inf}$} & \color{blue}{$\infty$} & \color{blue}{$\infty$} \\
  \hline 
  bi-Schr\"odinger ($s=2$) & 0 \footnotesize(under GCC) & \color{blue}{0} &  \color{blue}{$T_{\rm inf}$} & \color{blue}{$\infty$}\\
  \hline
  Heat& 0 & \color{blue}{$T_{\rm inf}$ or $\infty$} & ?& ?\\
  \hline
\end{tabular}
 \caption{Observability of subelliptic PDEs depending on the step.\newline If the results are established only in particular cases, they are in blue. \newline The first line is covered by Theorems \ref{t:obsRiemwaves} and \ref{t:main1}, the second line by Theorems \ref{t:lebschrod}, \ref{t:main2} and \ref{t:main3}, the third line by Theorem \ref{t:main2} and the fourth line by Theorems \ref{t:obsRiemheat}, \ref{t:bcg14} and \ref{t:koenig} (see also Corollary \ref{c:corheat}). For the last two interrogation marks, see Section \ref{s:conjheat}. Note that we illustrated Theorem \ref{t:main2} with the bi-Schr\"odinger equation, but we could have done it for a general $s$.}
  \label{t:summary}
\end{figure}
\end{center}

\subsection{Organization of the survey.} 
The goal of this survey is to give an overview of the ideas behind the three main results (and their proofs), namely Theorems \ref{t:main1}, \ref{t:main2} and \ref{t:main3}, to point out their common features, to develop ideas which were not necessarily written in the papers, and to explain how the tools developed along the proofs could be generalized. 

In Sections \ref{s:seul}, \ref{s:sun} and \ref{s:fermanian}, we explain respectively the main lines of the proofs of Theorems \ref{t:main1}, \ref{t:main2} and \ref{t:main3}. Therefore, Section \ref{s:seul} is quite geometric and presents for example the notion of nilpotentization of vector fields; Section \ref{s:sun} is more ``semi-classical'' and illustrates how resolvent estimates can be used to prove observability inequalities; and Section \ref{s:fermanian} is also ``semi-classical'', but the pseudo-differential operators used in this section are adapted to the group structure (that is, they come from non-commutative Fourier analysis). Sections \ref{s:sun} and \ref{s:fermanian}, although proving quite similar results, call for very different generalizations, which are discussed notably in Remark \ref{r:comparsuncfk}.

 In Section \ref{s:open}, we finally list a few natural open questions and directions of research which follow from our work.

\textbf{Acknowledgments.} 
I first thank Clotilde Fermanian Kammerer and Chenmin Sun for these very interesting collaborations, and my PhD advisors Emmanuel Tr\'elat and Yves Colin de Verdi\`ere who taught me so much about the subject. I also thank Clotilde Fermanian, Chenmin Sun and Emmanuel Tr\'elat for their suggestions and comments on a preliminary version of this survey. While working on the series of works reported here, I benefited from the interactions with many other mathematicians, among whom I am particularly grateful to Richard Lascar and V\'eronique Fischer for several interesting discussions. I was partially supported by the grant ANR-15-CE40-0018 of the ANR (project SRGI). 

\section{Subelliptic wave equations are never observable} \label{s:seul}
In this section, we explain ideas of proof for Theorem \ref{t:main1}. Let $T_0>0$ and $\omega\subset M$ be a measurable subset. We assume that there exist $1\leq i,j\leq m$ and $x$ in the interior of $M\backslash\omega$ such that $[X_i,X_j](x)\notin {\rm Span}(X_1(x),\ldots,X_m(x))$.  Under these assumptions, Theorem \ref{t:main1} will follow from Propositions \ref{p:existencegeod} and \ref{p:exactgb}. The first of these propositions says that given any open set of $M$ containing $x$, there exists a normal geodesic of the Hamiltonian attached to $X_1,\ldots,X_m$ traveled at speed $1$ which does not leave this open set on the time-interval $[0,T_0]$: this phenomenon is not true in Riemannian manifolds but is true in sub-Riemannian manifolds under the above assumptions. Thus, this normal geodesic remains far from $\omega$ on the time-interval $[0,T_0]$. The second proposition tells us that, as in the elliptic setting, it is possible to construct a sequence of solutions of the wave equation whose energy concentrates along this geodesic. This last fact contradicts the observability inequality.

Before stating these propositions, let us mention that, as in the elliptic setting, a normal geodesic is the projection of a null-bicharacteristic associated to the principal symbol $p_2$ of $\Delta$ (see Definition \ref{d:normal}).
\begin{prop} \label{p:existencegeod}
For any $T_0>0$, any $x\in M$ such that $[X_i,X_j](x)\notin {\rm Span}(X_1(x),\ldots,X_m(x))$ and any open neighborhood $V$ of $x$ in $M$ (with the initial topology on $M$), there exists a non-stationary normal geodesic $t\mapsto x(t)$ (traveled at speed $1$) such that $x(t)\in V$ for any $t\in [0,T_0]$.
\end{prop}
\begin{prop} \label{p:exactgb}
Let $[0,T_0]\ni t\mapsto x(t)$ be a non-stationary normal geodesic (traveled at speed $1$) which does not meet $\partial M$. 
Then there exists a sequence $(u_k)_{k\in\N^*}$ a sequence of solutions of the wave equation \eqref{e:system} such that
\begin{itemize}
\item The energy of $u_k$ is bounded below with respect to $k$ and $t\in [0,T_0]$:
\begin{equation} \label{e:convenergyforu2}
\exists A>0, \forall t\in [0,T_0],\quad  \liminf_{\substack{k\rightarrow +\infty}}E(u_k(t,\cdot))\geq A.
\end{equation}
\item The energy of $u_k$ is small off $x(t)$: for any $t\in [0,T_0]$, we fix $V_t$ an open subset of $M$ for the initial topology of $M$, containing $x(t)$, so that the mapping $t\mapsto V_t$ is continuous ($V_t$ is chosen sufficiently small so that this makes sense in a chart). Then
\begin{equation} \label{e:decreaseenergy2}
\sup_{t\in [0,T_0]}\int_{M\backslash V_t} \left(|\partial_tu_k(t,x)|^2+\sum_{j=1}^m (X_ju_k(t,x))^2\right)d\mu(x)\underset{k\rightarrow +\infty}{\rightarrow} 0.
\end{equation}
\end{itemize}
\end{prop}
The ``easy part'' is Proposition \ref{p:exactgb}, while the proof of Proposition \ref{p:existencegeod} is more complicated. Therefore, we start with a sketch of proof for Proposition \ref{p:exactgb} in Section \ref{s:exactgb}, and we give ideas for the proof of Proposition \ref{p:existencegeod} in Section \ref{s:existencegeod}.

\subsection{Gaussian beams in sub-Riemannian geometry}\label{s:exactgb}
We set $P=\partial_{tt}^2-\Delta$ and we consider the associated Hamiltonian
\begin{equation*}
p_2(t,\tau,x,\xi)=-\tau^2+g^*(x,\xi)
\end{equation*}
with $\tau$ the dual variable of $t$ and $g^*$ the Hamiltonian (or principal symbol) associated to $-\Delta$. For $\xi\in T^*M$, we have 
\begin{equation*}
g^*=\sum_{i=1}^m h_{X_i}^2.
\end{equation*}
Here, given any smooth vector field $X$ on $M$, we denoted by $h_X$ the Hamiltonian function (momentum map) on $T^*M$ associated with $X$ defined in local $(x,\xi)$-coordinates by $h_X(x,\xi)=\xi(X(x))$. Then $g^*$ is both the principal symbol of $-\Delta$, and also the cometric associated with $g$. 

In $T^*(\R\times M)$, the Hamiltonian vector field $\vec{H}_{p_2}$ associated with $p_2$ is given by $\vec{H}_{p_2}f=\{p_2,f\}$. Since $\vec{H}_{p_2}p_2=0$, we get that $p_2$ is constant along the integral curves of $\vec{H}_{p_2}$. Thus, the characteristic set $\mathcal{C}(p_2)=\{p_2=0\}$ is preserved under the flow of $\vec{H}_{p_2}$. Null-bicharacteristics are then defined as the maximal integral curves of $\vec{H}_{p_2}$ which live in $\mathcal{C}(p_2)$.  In other words, the null-bicharacteristics are the maximal solutions of
\begin{equation} \label{e:bicarac1}
\left\lbrace \begin{array}{l}
\dot{t}(s)=-2\tau(s)\,,\\
\dot{x}(s)=\nabla_\xi g^*(x(s),\xi(s))\,,\\
\dot{\tau}(s)=0\,,\\
\dot{\xi}(s)=-\nabla_xg^*(x(s),\xi(s))\,,\\
\tau^2(0)=g^*(x(0),\xi(0)).
\end{array}\right.
\end{equation}
This definition needs to be adapted when the null-bicharacteristic meets the boundary $\partial M$, but in the sequel, we only consider solutions of \eqref{e:bicarac1} on time intervals where $x(t)$ does not reach $\partial M$.

In the sequel, we take $\tau=-1/2$, which gives $g^*(x(s),\xi(s))=1/4$. This also implies that $t(s)=s+t_0$ and, taking $t$ as a time parameter, we are led to solve
\begin{equation} \label{e:bicarac}
\left\lbrace \begin{array}{l}
\dot{x}(t)=\nabla_\xi g^*(x(t),\xi(t))\,,\\
\dot{\xi}(t)=-\nabla_xg^*(x(t),\xi(t))\,,\\
g^*(x(0),\xi(0))=\frac{1}{4}.
\end{array}\right.
\end{equation}
In other words, the $t$-variable parametrizes null-bicharacteristics in a way that they are traveled at speed $1$.

\begin{definition} \label{d:normal}
A normal geodesic is the projection $t\mapsto x(t)$ of a null-bicharacteristic, (i.e., a solution of \eqref{e:bicarac}).
\end{definition}
Normal geodesics live in the ``elliptic part'' of $g^*$, i.e., where $g^*\neq 0$; this is the key point in the proof of Proposition \ref{p:exactgb}. Indeed, this result is well-known in the elliptic context, it is due to Ralston \cite{ralston1982gaussian}, and already H\"ormander noted that his argument extended to non-elliptic operators, as long as we were working in the elliptic part of the symbol (see \cite[Chapter 24.2]{hormander2007analysis}).

Let us start the proof of Proposition \ref{p:exactgb}. Taking charts of $M$, we can assume $M\subset\R^n$. In the sequel, we change a bit the notations: we use $x=(x_0,x_1,\ldots,x_n)$ where $x_0=t$ in the earlier notations, and we set $x'=(x_1,\ldots,x_n)$. Similarly, we take the notation $\xi=(\xi_0,\xi_1,\ldots,\xi_n)$ where $\xi_0=\tau$ previously, and $\xi'=(\xi_1,\ldots,\xi_n)$. The bicharacteristics are parametrized by $s$ as in \eqref{e:bicarac1}, and without loss of generality, we only consider bicharacteristics with $x(0)=0$ at $s=0$, which implies in particular $x_0(s)=s$ because of our choice $\tau^2(s)=g^*(x(s),\xi(s))=1/4$. In the sequel, a null-bicharacteristic $s\mapsto (x(s),\xi(s))$ is fixed, with $x(0)=0$.
 
We take
\begin{equation} \label{e:ansatz}
v_k(x)=k^{\frac{n}{4}-1}a_0(x)e^{ik\psi(x)}.
\end{equation}
where the phase $\psi(x)$ is quadratic,
\begin{equation*}
\psi(x)=\xi'(s)\cdot (x'-x'(s))+\frac12 (x'-x'(s))\cdot M(s)(x'-x'(s)),
\end{equation*} 
for $x=(t,x')\in\R\times\R^{n}$ and $s$ such that $t=t(s)$. This choice of $v_k$ is the so-called WKB ansatz, and it only yields \emph{approximate solutions} of the wave equation \eqref{e:waveop}. Indeed, there holds
\begin{equation} \label{e:opondesvk}
\partial_{tt}^2v_k-\Delta v_k=(k^{\frac{n}{4}+1}A_1+k^{\frac{n}{4}}A_2+k^{\frac{n}{4}-1}A_3)e^{ik\psi}
\end{equation}
with 
\begin{align}
A_1(x)&=p_2\left(x,\nabla\psi(x)\right)a_0(x) \nonumber \\
A_2(x)&=La_0(x)\label{e:A2} \qquad \text{($L$ is a linear transport operator)}\nonumber\\
A_3(x)&=\partial_{tt}^2a_0(x)-\Delta a_0(x).\nonumber
\end{align}
If we take $\psi$ to be complex-valued (by choosing a complex-valued matrix $M(s)$), then $e^{ik\psi}$ looks like a Gaussian centered at $x'(s)$, for any $s$. Hence, for large $k$, we see on \eqref{e:opondesvk} is naturally small outside the geodesic $s\mapsto x'(s)$, and the only place where it may be large is precisely on the geodesic. In order for $v_k$ to be an approximate solution of the wave equation, it is thus sufficient to have $A_1, A_2, A_3$ vanish at sufficiently high order on the geodesic curve, and this is achieved by choosing adequately $M(s)$ (which is a complex-valued matrix varying continuously with $s$) and $a_0(x)$.

More precisely, in order to achieve the quantitative bound $\|\partial_{tt}^2v_k-\Delta v_k\|_{L^1((0,T);L^2(M))}\leq Ck^{-\frac12}$, one has to 
\begin{itemize}
\item take $M$ as a solution of the Riccati equation
\begin{equation*}
\frac{dM}{ds}+MCM+B^TM+MB+A=0
\end{equation*}
where $A,B,C$ are explicit functions of $s$ which can be expressed in terms of the second derivatives of $p_2$ along $s\mapsto (x(s),\xi(s))$;
\item take $a_0$ so that $a_0(x(0))\neq 0$ and $La_0$ vanishes along $s\mapsto x(s)$, i.e., $La_0(x(s))=0$. This is possible since $L$ is a linear transport operator.
\end{itemize}
The first point yields that $A_1$ together with its first and second derivatives vanish along the geodesic curve $s\mapsto x(s)$, and the second point implies that $A_2$ vanishes along this same curve. Beside implying that $v_k$ is an approximate solution of the wave equation, these choices guarantee that the energy of $v_k$ concentrates (uniformly in $s$) along the geodesic:
\begin{equation*}
\sup_{t\in [0,T]}\int_{M\backslash V_t} \left(|\partial_tv_k(t,x)|^2+\sum_{j=1}^m (X_jv_k(t,x))^2\right)d\mu(x)\underset{k\rightarrow +\infty}{\rightarrow} 0
\end{equation*}
where $V_t$ is defined in Proposition \ref{p:exactgb}. Note also that the bound
\begin{equation*}
\exists A>0, \forall t\in [0,T_0],\quad  \liminf_{\substack{k\rightarrow +\infty}}E(v_k(t,\cdot))\geq A
\end{equation*}
is satisfied thanks to the normalizing constant $k^{\frac{n}{4}-1}$ in the definition of $v_k$ \eqref{e:ansatz}.

In order to pass from approximate to exact solutions of the wave equation, one chooses $u_k$ to have the same initial data as $v_k$ and to be an exact solution of the wave equation \eqref{e:waveop}. Then, the Gronwall lemma ensures that \eqref{e:convenergyforu2} and \eqref{e:decreaseenergy2} are satisfied, and Proposition \ref{p:exactgb} is proved.

Let us end this section with two other possible points of view on Proposition \ref{p:exactgb} and its above proof. First, it can be reformulated in terms of propagation of Lagrangian spaces, as written in \cite{letrouit2020subelliptic}. This point of view was developed for example by H\"ormander in \cite[Chapter 24]{hormander2007analysis}. Another reformulation is that microlocal defect measures propagate in the elliptic part of the symbol as for the usual elliptic wave equation: in other words, one could take a sequence of initial data concentrating microlocally on the starting point $(x(0),\xi(0))$ of the geodesic and prove that the microlocal defect measure associated to the solutions propagates following the geodesic flow (see \cite{gerard1991microlocal}). But the Gaussian beam construction is interesting in its own and particularly simple, this is why we presented this point of view here.

\subsection{Spiraling geodesics} \label{s:existencegeod} Proposition \ref{p:existencegeod} can be easily seen to hold in the Heisenberg group, as shown in Section \ref{s:heiscase}; the proof then consists in extending its validity to larger classes of sub-Laplacians, until reaching the level of generality of Theorem \ref{t:main1}. 
\subsubsection{Spiraling geodesics in the 3D flat Heisenberg case.} \label{s:heiscase}
We consider the three-dimensional manifold with boundary $M_{1}=(-1,1)_{x_1}\times \mathbb{T}_{x_2}\times \mathbb{T}_{x_3}$, where $\mathbb{T}=\R/\mathbb{Z}\approx (-1,1)$ is the 1D torus. We endow $M_{1}$ with the vector fields $X_1=\partial_{x_1}$ and $X_2=\partial_{x_2}-x_1\partial_{x_3}$ and we consider the associated sub-Laplacian $\Delta=X_1^2+X_2^2$. This sub-Riemannian structure is called the ``Heisenberg manifold with boundary", it is a variant of Example \ref{exHeis} which has no boundary. We endow $M$ with an arbitrary smooth volume $\mu$ and we denote by 
\begin{equation} \label{e:heishamilt}
g^*_{\text{Heis}}=\xi_1^2+(\xi_2-x_1\xi_3)^2
\end{equation}
the 3D flat Heisenberg Hamiltonian.

The geodesics we consider are given by
\begin{equation*}
\begin{array}{l}
x_1(t)=\varepsilon\sin(t/\varepsilon)\\
x_2(t)=\varepsilon\cos(t/\varepsilon)-\e \\
x_3(t)=\varepsilon(t/2-\varepsilon\sin(2t/\varepsilon)/4).
\end{array}
\end{equation*}
They spiral around the $x_3$ axis $x_1=x_2=0$.

Here, one should think of $\varepsilon$ as a small parameter. In the sequel, we denote by $x_\varepsilon$ the geodesic with parameter $\varepsilon$. The associated momenta are 
\begin{equation} \label{e:momenta}
\xi_{1}=\frac12\cos(t/\e),\quad  \xi_{2}=0 \quad \text{and} \quad \xi_{3}=\frac{1}{2\varepsilon},
\end{equation}
 and we can check that that $g^*_{\text{Heis}}\equiv 1/4$. The constant $\xi_{3}$ is a kind of rounding number reflecting the fact that the geodesic spirals at a certain speed around the $x_3$ axis. To obtain a geodesic which makes smaller spirals, we choose a larger covector very close to $\mathcal{D}^\perp$.

To prove Proposition \ref{p:existencegeod} in the case of the Heisenberg manifold, we can assume without loss of generality that $V$ contains $0$. Then, given any $T_0>0$, for $\varepsilon$ sufficiently small, we have $x_\varepsilon(t)\in V$ for every $t\in (0,T_0)$. This proves Proposition \ref{p:existencegeod} in this case.

\subsubsection{Spiraling when length $\geq 3$ brackets vanish} \label{s:bracketslength3}
We now explain how to prove Proposition \ref{p:existencegeod} in a slightly more general case: in this paragraph, we assume that $[X_i,[X_j,X_k]]=0$ for any $1\leq i,j,k\leq m$. 

More precisely, we assume that $M\subset \R^n$ (with coordinates $x_1,\ldots,x_n$), and that for any $1\leq j\leq m$, 
\begin{equation} \label{e:reducedvf}
X_j=\sum_{j=1}^n a_{ij}\partial_{x_i}
\end{equation}
where $a_{ij}$ is a constant when $i\leq m$, and $a_{ij}=cx_\ell+d$ when $i\geq m+1$, for some $\ell\leq m$ that may depend on $i$ and $j$. One can verify that $[X_i,[X_j,X_k]]=0$ for any $1\leq i,j,k\leq m$.

Our goal is to isolate a direction which will play the role of the direction $\xi_3$ of large covectors in \eqref{e:momenta}. A similar spiraling as in the above Heisenberg case happens when we take covectors in an invariant plane of the Goh matrix (thus the associated control describes circles in this invariant plane). Let us explain this phenomenon in detail.

In its ``control'' form, the equation of normal geodesics can be written as follows (we refer the reader to \cite[Chapter 4]{agrachev2019comprehensive}):
\begin{equation}\label{e:controlform}
\dot{x}(t)=\sum_{i=1}^m u_i(t)X_i(x(t)),
\end{equation}
where the $u_i$ are the controls, explicitly given by
\begin{equation}\label{e:uexplicit}
u_i(t)=2h_{X_i}(x(t),\xi(t)).
\end{equation} 
Thanks to \eqref{e:reducedvf}, we rewrite \eqref{e:controlform} as
\begin{equation}\label{e:diffeqxt}
\dot{x}(t)=F(x(t))u(t), 
\end{equation}
where $F=(a_{ij})$, which has size $n\times m$, and $u=\;^t(u_1,\ldots,u_m)$. Differentiating \eqref{e:uexplicit}, we have the complementary equation
$$
\dot{u}(t)=G(x(t),\xi(t))u(t)
$$
where $G$ is the Goh matrix
$$
G=(2\{h_{X_j},h_{X_i}\})_{1\leq i,j\leq m}
$$
(it differs from the usual Gox matrix by a factor $-2$ due to the absence of factor $\frac12$ in the Hamiltonian $g^*$ in our notations). One can check thanks to \eqref{e:reducedvf} that $G(t)$ is constant in $t$. 

\smallskip

We know that $G\neq 0$ and that $G$ is antisymmetric. The whole control space $\R^m$ is the direct sum of the image of $G$ and the kernel of $G$, and $G$ is nondegenerate on its image.   We take $u_0$ in an invariant plane of $G$; in other words its projection on the kernel of $G$ vanishes (see Remark \ref{r:singular}). We denote by $\widetilde{G}$ the restriction of $G$ to this invariant plane. We also assume that $u_0$, decomposed as $u_0=(u_{01},\ldots,u_{0m})\in\R^m$, satisfies $\sum_{i=1}^m u_{0i}^2=1/4$. Then $u(t)=e^{t\widetilde{G}}u_0$ and since $e^{t\widetilde{G}}$ is an orthogonal matrix, we have $\|e^{t\widetilde{G}}u_0\|=\|u_0\|$. We have by integration by parts
\begin{align}
x(t)&=\int_0^tF(x(s))e^{s\widetilde{G}}u_0\,ds\nonumber\\
&=F(x(t))\widetilde{G}^{-1}(e^{t\widetilde{G}}-I)u_0-\int_0^t\frac{d}{ds}(F(x(s))\widetilde{G}^{-1}(e^{s\widetilde{G}}-I)u_0\,ds.\label{e:duhamel}
\end{align}

\smallskip

Let us now choose the initial data of our family of normal geodesics (indexed by $\e$). The starting point $x^\e(0)=0$ is the same for any $\e$, we only have to specify the initial covectors $\xi^\e=\xi^\e(0)\in T_0^*\R^n$. For any $i=1,\ldots,m$, we impose that
\begin{equation}\label{e:covyieldscontrol}
\langle \xi^\e, X_i\rangle=u_{0i}.
\end{equation}
It follows that $g^*(x(0),\xi^\e(0))=\sum_{i=1}^m u_{0i}^2=1/4$ for any $\e>0$. Now, we notice that $\text{Span}(X_1,\ldots,X_m)$ is in direct sum with the Span of the $[X_i,X_j]$ for $i,j$ running over $1,\ldots,m$ (this follows from \eqref{e:reducedvf}). Fixing $G^0\neq 0$ an antisymmetric matrix and $\widetilde{G}^0$ its restriction to an invariant plane, we can specify, simultaneously to \eqref{e:covyieldscontrol}, that
$$
\langle \xi^\e, 2[X_j,X_i]\rangle = \varepsilon^{-1}G^0_{ij}.
$$
Then $x^\e(t)$ is given by \eqref{e:duhamel} applied with $\widetilde{G}=\e^{-1}\widetilde{G}^0$, which brings a factor $\e$ in front of \eqref{e:duhamel}.

\smallskip

Recall finally that the coefficients $a_{ij}$ which compose $F$ are degree $1$ (or constant) homogeneous polynomials in $x_1,\ldots,x_{m}$. Thus $\frac{d}{ds}(F(x(s))$ is a linear combination of $\dot{x}_i(s)$ which we can rewrite thanks to \eqref{e:diffeqxt} as a combination with bounded coefficients (since $\sum_{i=1}^m u_i^2=1/4$) of the $x_i(s)$. Hence, applying the Gronwall lemma in \eqref{e:duhamel}, we get $\|x^\e(t)\|\leq C\varepsilon$, which concludes the proof of Proposition \ref{p:existencegeod} in this case.
\begin{remark} \label{r:singular}
Let us explain why we choose $u_0$ to be in an invariant plane of $G$. If the projection of $u_0$ to the kernel of $G$ is nonzero then
the primitive of the exponential of $e^{\frac t\varepsilon G_0}u_0$ contains a linear term that does not depend on $\varepsilon$. Then the corresponding trajectory follows a singular curve (see \cite[Chapter 4]{agrachev2019comprehensive} for a definition). This means, we find normal geodesics which spiral around a singular curve and do not remain close to their initial point over $(0,T_0)$, although their initial covector is ``high in the cylinder bundle $U^*M$''. For example, for the Hamiltonian $\xi_1^2+(\xi_2+x_1^2\xi_3)^2$  associated to the ``Martinet'' vector fields $X_1=\partial_{x_1}$, $X_2=\partial_{x_2}+x_1^2\partial_{x_3}$ in $\R^3$, there exist normal geodesics which spiral around the singular curve $(t,0,0)$.
\end{remark}

\subsubsection{Spiraling in the general case} \label{s:nilpo}
The reduction of the general case to the case of Section \ref{s:bracketslength3} where all length $\geq 3$ brackets vanish is done through the nilpotentization procedure which dates back to \cite{rothschild1976hypoelliptic}.  The reader can refer to \cite[Chapter 10]{agrachev2019comprehensive} and \cite[Chapter 2]{jean2014control} for recent introductions to the subject.


Essentially, the nilpotentization procedure consists in a truncation in the Taylor series of the vector fields $X_1,\ldots,X_m$ which define the sub-Laplacian.  We will not describe the nilpotentization procedure in details here, but just give an example.

 \begin{example} \label{e:nilpo} We reproduce here the example \cite[Example 2.8]{jean2014control}. We consider the vector fields on $\R^2\times\T$ given by $X_1=\cos(\theta)\partial_x+\sin(\theta)\partial_y$ and $X_2=\partial_\theta$. We have $[X_1,X_2]=\sin(\theta)\partial_x-\cos(\theta)\partial_y$. At $q=0$, we have $X_1(0)=\partial_x$, $X_2(0)=\partial_\theta$ and $[X_1,X_2](0)=-\partial_y$. Thus, we say that $x$ and $\theta$ have ``weight $1$'', while $y$ has ``weight $2$'', because the coordinate $y$ ``needs a bracket to be generated''. And we also attribute weights to vector fields: $\partial_x$ and $\partial_\theta$ have weight $-1$ while $\partial_y$ has weight $-2$. The rule is that the weight of a product is the sum of the weights: for example, $x\partial_\theta$ has weight $1+(-1)=0$.
 
Now we write the Taylor expansion of $X_1$ and $X_2$ in the coordinates $(x,\theta,y)$ in the form $X_1=X_1^{(-1)}+X_1^{(0)}+X_1^{(1)}+\ldots$ where each $X_1^{(k)}$ has weight $k$ (and similarly for $X_2$). This Taylor expansion of $X_1$ and $X_2$ at $q=0$ yields the homogeneous components
\begin{equation*}
X_1^{(-1)}=\partial_x+\theta\partial_y, \qquad X_1^{(0)}=0, \qquad X_1^{(1)}=-\frac{\theta^2}{2}\partial_x-\frac{\theta^3}{6}\partial_y, \qquad \ldots
\end{equation*}
and $X_2^{(-1)}=X_2=\partial_\theta$. We define the nilpotent approximation of $(X_1,X_2)$ at $q=0$ in coordinates $(x,\theta,y)$ to be the vector fields $\widehat{X}_1$ and $\widehat{X}_2$ given by the ``main order terms'', namely 
\begin{equation*}
\widehat{X}_1=X_1^{(-1)}=\partial_x+\theta\partial_y, \qquad \widehat{X}_2=X_2^{(-1)}=\partial_\theta.
\end{equation*}
These vector fields generate a nilpotent Lie algebra of step $2$ since all brackets of length $\geq 3$ between $\widehat{X}_1$ and $\widehat{X}_2$ vanish.
\end{example}

Let us explain how finish the proof of Proposition \ref{p:existencegeod} in the general case. We fix $q\in M$. Thanks to Section \ref{s:bracketslength3} we can find a normal geodesic associated to the nilpotentized (at $q$) Hamiltonian
\begin{equation} \label{e:nilpohamilt}
\widehat{g}^*=\sum_{i=1}^mh_{\widehat{X}_i}^2
\end{equation}
which stays very close to $q$. Then, the key argument is that since we work locally near $q$ and since the vector field $\widehat{X}_i$ is a very good approximation of $X_i$ near $q$ for any $i$, the normal geodesics associated to the initial Hamiltonian $g^*$ cannot be far from those of $\widehat{g}^*$ when working near $q$. In other words, there exists a normal geodesic associated to the initial Hamiltonian 
$g^*$ which stays very close to $q$. 

One should however take care that the situation is not always as favorable as in Example \ref{e:nilpo}, for two reasons:
\begin{itemize}
\item To obtain good properties of the ``truncated'' vector fields, one needs to do the Taylor expansion in a good system of coordinates. In the above example, it was quite ``natural'' to take $(x,\theta,y)$ coordinates, but in general one should work in so-called ``privileged coordinates'' to guarantee that the geodesics of the nilpotentized Hamiltonian are not too far from those of the initial Hamiltonian.
\item The nilpotentized vector fields always form a nilpotent system, meaning that there exists $k\in\N$ such that all length $\geq k$ brackets between $\widehat{X}_1,\ldots,\widehat{X}_m$ vanish. But in general, $k$ is not necessarily equal to $3$ (as was assumed in Section \ref{s:bracketslength3}). To reduce to the case $k=3$, one has again to ``compare'' the geodesics of the nilpotentized Hamiltonian with a simpler Hamiltonian, defined with a system of vector fields such that all their length $\geq 3$ brackets vanish.
\end{itemize}
To sum up, the proof of Proposition \ref{p:existencegeod} goes by successive simplifications of the vector fields, until we arrive at the situation of Section \ref{s:bracketslength3} which can be handled ``by hand'' (i.e., we can find explicitly the initial covectors of the geodesics).

\section{Subellipticity and strength of propagation: Baouendi-Grushin Schr\"odinger equations} \label{s:sun}
\subsection{Motivation.} It is well-known that the Riemannian wave equation propagates at finite speed and that the Riemannian Schr\"odinger equation propagates at infinite speed. For subelliptic equations, the propagation ``at null speed'' of the wave equation (shown to be a general phenomenon in Section \ref{s:seul}), and the existence of travelling waves solutions of the Heisenberg Schr\"odinger equation (see \cite{bahouri2000espaces}) motivated us to undertake a general study of the propagation speed of subelliptic equations depending on two parameters: the step $k$ and the strength of propagation $s$. This last parameter, which will be defined below, is equal to $1/2$ for the wave equation and to $1$ for the Schr\"odinger equation. 

As a first step in this study, we took the model family of operators given by Example \ref{exgamma}, i.e., we consider for $\gamma \geq 0$ (not necessarily an integer) the sub-Laplacian $\Delta_\gamma=\partial_x^2+|x|^{2\gamma}\partial_y^2$ on the manifold $M=(-1,1)_x\times\T_y$. In case $\gamma\in\N$, the step is given by $k=\gamma+1$. Then, for $s\in\N$, we consider the equation \eqref{e:schrodfrac}. The strength of propagation is the parameter $s$ appearing in this equation; this terminology seemed natural to us since equations with large $s$ tend to propagate more quickly. The domain of observation is a strip $\omega=(-1,1)\times I$ where $I\subset \T$ as mentioned in the introduction. Depending on $\gamma$ and $s$, Theorem \ref{t:main2} says if the observability inequality \eqref{e:obsschrod} holds or not: in particular, it makes appear a threshold $s=\frac{\gamma+1}{2}$ at which subellipticity effects and propagation effects balance each other.

The proof of Theorem \ref{t:main2} splits into two parts:
\begin{enumerate}
\item To disprove observability when $s<\frac12(\gamma+1)$ or when $s=\frac12(\gamma+1)$ and time is small (i.e., for Point (3) and part of Point (2)), we construct solutions of \eqref{e:schrodfrac} which propagate along the vertical line $\{x=0\}$. This is different from what was done in Section \ref{s:seul}: in Section \ref{s:seul}, we were working in the elliptic part of the symbol of the sub-Laplacian, whereas here the constructed solutions propagate in the characteristic manifold (where the sub-Laplacian is not elliptic). Thus, the present construction is more involved and less general than that of Section \ref{s:seul}, but these difficulties are unavoidable since we are dealing with Schr\"odinger equations whose speed of propagation in the elliptic part of the symbol is infinite. This Gaussian beam construction is described in Section \ref{s:vgb}
\item To prove observability when $s>\frac12(\gamma+1)$ or when $s=\frac12(\gamma+1)$ and time is large (i.e., for Point (1) and part of Point (2)), we establish a sharp resolvent estimate, i.e., a time-independent inequality which describes the maximal concentration outside $\omega$ of an approximate eigenfunction of $\Delta_\gamma$. By a classical argument due to Nicolas Burq and Maciej Zworski (see \cite{burq2004geometric}), this implies observability. Indeed, this method of proof gives as corollaries an observability result for the heat equation and a decay estimate for the damped wave equation, both driven by the sub-Laplacian $\Delta_\gamma$.  The resolvent estimate is presented and commented in Section \ref{s:resgrushin}, and (sketchly) proved in Section \ref{s:proofres}.
\end{enumerate}

\subsection{Vertical Gaussian beams for contradicting observability.} \label{s:vgb} 

For proving Point (3) and Point (2) for small times, it is indeed sufficient to deal with the ``critical'' case $s=\frac12(\gamma+1)$. Point (3) then follows immediately from the abstract result \cite[Corollary 3.9]{miller2012resolvent}: if \eqref{e:schrodfrac} was observable for some $T>0$ and some $s<\frac{\gamma+1}{2}$, then it would be observable in any time for $s=\frac12(\gamma+1)$, which would be a contradiction. Hence, in the sequel, we assume $s=\frac12(\gamma+1)$. Moreover, we also assume that $\gamma>1$; for the case $s=\gamma=1$, the argument is slightly different, see \cite{burq2019time}.  

The non-observability part of Point (2) immediately follows from the following proposition:
\begin{prop} \label{p:gbpoint3}
There exist $T_0>0$ and a sequence of solutions $(v_n)_{n\in\N}$ of \eqref{e:schrodfrac} with initial data $(v_n^0)_{n\in\N}$  such that $\|v_n^0\|_{L^2(M)}=1$ and
\begin{equation} \label{e:smallmass}
\int_0^{T_0}\int_\omega |v_n(t,x,y)|^2dxdydt \underset{n\rightarrow +\infty}{\longrightarrow} 0.
\end{equation}
\end{prop}
Let us explain the intuition behind this result on the simpler example of the Grushin-Schr\"odinger equation
\begin{equation} \label{e:Grushschrod}
i\partial_tu+\Delta_G u=0
\end{equation}
where $\Delta_G=\partial_x^2+x^2\partial_y^2$ has been introduced in Example \ref{exGrushin}. We consider this equation in $\R_t\times\R_x\times\R_y$ instead of $\R_t\times (-1,1)_x\times\T_y$ since computations become simpler in this (non-compact) setting. Making a Fourier expansion in the $y$ variable, we see that we obtain from $-\Delta_G$ a family of harmonic oscillators $-\partial_x^2+x^2\eta^2$ whose associated eigenvalues are $(2m+1)|\eta|$, $m\in\N$. In other words, 
\begin{equation} \label{e:decompoL2grushin}
L^2(\R^2)=\underset{\pm}{\oplus}\underset{m\in\N}{\oplus} V_m^\pm, \qquad \Delta_{G|V_m^\pm}=\pm i(2m+1)\partial_y,
\end{equation}
and thus the solution of \eqref{e:Grushschrod} is obtained by solving an infinite number of transport equations along the $y$-axis, traveling at speed $2m+1$ for $m\in\N$. More explicitly, the solution of \eqref{e:Grushschrod} in  $\R_t\times\R_x\times\R_y$ is
\begin{equation*}
u(t,x,y)=\frac{1}{\sqrt{2\pi}}\sum_{m=0}^\infty \int_{\R} e^{-it(2m+1)|\eta|+iy\eta}\widehat{u}_{0,m}(\eta)h_m(\sqrt{|\eta|}x)d\eta
\end{equation*}
where $h_m$ is the $m$-th eigenfunction of the harmonic oscillator $-\partial_x^2+x^2$ on $\R$ (Hermite function). Our goal is to find solutions whose energy is concentrated in $M\setminus \omega$ during a long time, therefore these solutions should travel slowly along the $y$-axis. It will be the case if we manage to construct solutions corresponding to low values of $m$; indeed Burq and Sun \cite[Section 9]{burq2019time} showed how to construct such solutions whose only non-vanishing component is the ``mode'' $m=0$ (but $\eta$ takes different values).

The case of the equation $i\partial_tu-(-\Delta_\gamma)^su=0$ (or, more precisely, \eqref{e:schrodfrac}) is more involved but it is based on the same idea. The computations are less explicit since the eigenfunctions of $\Delta_\gamma$ are less explicit than those of $\Delta_G$, but the knowledge of the behaviour at infinity of the first eigenfunction of the harmonic oscillator is indeed sufficient to do the computations. This (positive) normalized ground state of the operator $Q_\gamma=-\partial_z^2+|z|^{2\gamma}$ on $\R_z$ is denoted by $\phi_\gamma$ and it satisfies
\begin{equation*}
Q_\gamma\phi_\gamma=\mu_0\phi_\gamma.
\end{equation*}
The normalized ground state of the operator $P_{\gamma,w}=-\partial_x^2+|x|^{2\gamma}w^2$ on $\R_x$ is then
\begin{equation*}
p_\gamma(w,x)=|w|^{\frac{1}{2(\gamma+1)}}\phi_\gamma(|w|^{\frac{1}{\gamma+1}}x).
\end{equation*}
and the associated eigenvalue is $\lambda_\gamma(w)=\mu_0|w|^{\frac{2}{\gamma+1}}$. The relation between the different variables is $z=|w|^{\frac{1}{\gamma+1}}x$, and $w$ will be taken to equal to $\eta$, the Fourier dual variable of $y$.

 We assume without loss of generality that $\omega=\R_x\times I_y$, where $I=(a,b)$, with $0<a<b\leq +\infty$. Let us fix $T_0<a/\mu_0^s$.

We take $\psi\in C_c^{\infty}(\frac{1}{2}\leq |\eta|\leq 1)$, a sequence $h_n\rightarrow 0$, and we consider
\begin{equation*}
 v_n(t,x,y)=h_n^{\frac12}\int_{\R}\psi(h_n\eta)\;e^{iy\eta-it\mu_0^s|\eta|}p_\gamma(\eta,x)\;d\eta,
\end{equation*}
which is a solution of \eqref{e:schrodfrac} (recall that $s=(\gamma+1)/2$). It has only high $\eta$-modes thanks to the cutoff $\psi(h_n\eta)$, and each of these modes travels at speed $\mu_0^s$ along the $y$-axis since $2s/(\gamma+1)=1$. By the Plancherel theorem, we can see that $\|v_{n,0}\|_{L^2}\gtrsim 1$, so that we just have to prove \eqref{e:smallmass}.

The key intuition is that, for any $t\in\R_+$, the mass of $v_n(t,\cdot,\cdot)$ is concentrated near $(x,y)=(0,\mu_0^st)$, and thus outside $\omega$ for $t\leq T_0$. To prove it, one uses the  Poisson formula, which yields
\begin{equation*}
v_n(t,x,y)=\sum_{m\in\mathbb{Z}}  \widehat{K_{t,x,y}^{(n)}}(2\pi m)
\end{equation*}
where $\Phi_m(t,y,w)=yw-2\pi mw- t\mu_0^sw$ and
\begin{equation*}
 \widehat{K_{t,x,y}^{(n)}}(2\pi m)=h_n^{\frac12}\int_{\R}\psi(h_nw)p_\gamma(w,x)e^{i\Phi_m(t,y,w)}dw.
\end{equation*}
Then, \eqref{e:smallmass} follows from the fact that for $|y|\geq a$ and $t\leq T_0$, each $ \widehat{K_{t,x,y}^{(n)}}(2\pi m)$ is small.

\subsection{Resolvent estimate.}  \label{s:resgrushin} This subsection is devoted to explaining how time-independent estimates can be relevant for proving observability inequalities (which inherently involve a time variable), and to apply it to our particular goal of proving Point (1) and part of Point (2) of Theorem \ref{t:main2}. The starting point of the analysis is the following theorem due to Nicolas Burq and Maciej Zworski, which will be commented after its statement:
\begin{thm}[\cite{burq2004geometric}]\label{abstract}
	Let $P(h)$ be self-adjoint on some Hilbert space $\mathcal{H}$ with densely defined domain $\mathcal{D}$ and $A(h):\mathcal{D}\rightarrow \mathcal{H}$ be bounded. Assume that uniformly for $\tau\in I=[-b,-a]\subset \R$, we have the following resolvent inequality
\begin{equation}  \label{e:quasimodeest}
\|u\|_{\mathcal{H}}\leq \frac{G(h)}{h}\|(P(h)+\tau)u\|_{\mathcal{H}}+g(h)\|A(h)u\|_{\mathcal{H}}
\end{equation}
	for some $1\leq G(h)\leq O(h^{-N_0})$.
	Then there exist constants $C_0,c_0,h_0>0$, such that for every $T(h)$ satisfying
	$$ \frac{G(h)}{T(h)}<c_0,
	$$
	we have, for all $0<h<h_0$
	$$ \|\psi(P(h))u\|_{\mathcal{H}}^2\leq C_0\frac{g(h)^2}{T(h)}\int_0^{T(h)}\|A(h)e^{-\frac{itP(h)}{h}}\psi(P(h))u\|_{\mathcal{H}}^2dt,
	$$
	where $\psi\in C_c^{\infty}((a,b))$.
\end{thm}
This statement, proved in \cite[Section 3]{burq2004geometric}, calls for several comments:
\begin{itemize}
\item Typically, $P(h)=-h^2\Delta$, or, in our case $P(h)=-h^2\Delta_\gamma$ with domain
\begin{equation*}
D(\Delta_\gamma)=\{u\in \mathcal{D}'(M) : \partial_x^2 u, |x|^{2\gamma}\partial_y^2u\in L^2(M) \text{ and } u_{|\partial M}=0\}.
\end{equation*}
 Therefore, the term $(P(h)+\tau)u$ gives an estimate on how far from a true eigenfunction $u$ is. The small parameter $h$ is called semi-classical parameter, and it naturally shows up  when studying Schr\"odinger equations. 
\item The operator $A(h)$ is the observation operator, and in our case it does not depend on $h$: it is equal to the characteristic function $\mathbf{1}_\omega$ of $\omega$.
\item The expression $\psi(P(h))$ is defined thanks to functional calculus, and $T(h)$ has to be understood as the time-scale at which the observability inequality is valid. The parameter $\tau$ shows up naturally when taking the Fourier transform in time of the equation: it is the dual variable of $t$.
\end{itemize}

In our case, the resolvent estimate (or quasimode estimate) \eqref{e:quasimodeest} takes the following form:
\begin{thm} \label{t:resgamma}
Let $\gamma\geq 1$ and let $\omega$ contain a horizontal strip $(-1,1)\times I$.  There exist $C, h_0>0$ such that for any $v\in D(\Delta_\gamma)$ and any $0<h\leq h_0$, there holds 
\begin{equation*}
 \|v\|_{L^2(M)}\leq C(h^{-(\gamma+1)}\|(h^2\Delta_\gamma+1)v\|_{L^2(M)}+\|v\|_{L^2(\omega)}).
\end{equation*}
\end{thm}
The exponent $h^{-(\gamma+1)}$ is optimal. It may seem strange that the parameter $s$ does not appear in the above resolvent inequality. However, we can deduce directly from Theorem \ref{t:resgamma} another resolvent inequality where $s$ appears:
$$ \|u\|_{L^2(M)}\leq C(h^{-(\gamma+1)}\| ((-h^2\Delta_\gamma)^s-1)u\|_{L^2(M)}+\|u\|_{L^2(\omega)}).
$$
From this and Theorem \ref{abstract}, one obtains the following \emph{spectrally localized observability inequality}:
\begin{equation} \label{e:locobs}
 \|u_h\|_{L^2(M)}^2\leq C\int_0^{T}\|e^{-it'(-\Delta_\gamma)^s}u_h\|_{L^2(\omega)}^2dt'
\end{equation}
where $u_h=\psi((-h^2\Delta_\gamma)^s )u$ with $\psi\in C_c^\infty((1/2,2))$. This is not totally immediate, the precise argument is written in \cite[Section 3.1]{letrouit2020observability}. From \eqref{e:locobs}, we can conclude thanks to a standard unique continuation argument the proof of Theorem \ref{t:main2}, that is, we transform this inequality for $u_h$ (which is spectrally localized) to the inequality \eqref{e:obsschrod} which holds for any initial datum $u_0$. This unique continuation argument is based on the fact that if $u_*$ is an eigenfunction of $(-\Delta_\gamma)^s$ which vanishes on $\omega$, then it vanishes on the whole manifold $M$. We refer to \cite[Section 3.2]{letrouit2020observability} for full details.

\subsection{Proof of Theorem \ref{t:resgamma}.} \label{s:proofres} The proof of Theorem \ref{t:resgamma} is quite long, and we will only describe here its main lines. It consists in decomposing solutions according to the joint eigenspaces of $|D_y|$ and $-\Delta_\gamma$ (it has already been done in Section \ref{s:vgb}); then we distinguish between several ``regimes'' depending on the value of $|D_y|$ compared to the value of $-\Delta_\gamma$, and this corresponds to different geometric propagations of solutions. This strategy is inspired by that of \cite{burq2019time}, but here we work in a time-independent setting since we seek to prove a resolvent estimate. This fine analysis of the different microlocal regimes should be compared with the non-commutative approach of Section \ref{s:fermanian}, see Remark \ref{r:comparsuncfk}.

The proof is by contradiction.  Assume that there exists a sequence $(v_h)_{h>0}$ such that
\begin{equation} \label{e:seekcontrad}
\|v_h\|_{L^2(M)}=1, \quad \|v_h\|_{L^2(\omega)}=o(1), \quad \|f_h\|_{L^2(M)}=o(h^{\gamma+1})
\end{equation}
where $f_h=(h^2\Delta_\gamma+1)v_h$, and we seek for a contradiction, which would prove Theorem \ref{t:resgamma}. A short argument shows that we can furthermore assume that $v_h=\psi(h^2\Delta_\gamma)v_h$ where $\psi\in C_c^\infty((-\infty,0))$ verifies: $\psi\equiv 1$ near $-1$ and $\psi=0$ outside $(-2,-\frac{1}{2})$. Here and in all the sequel, we use functional calculus to define expressions such as $\psi(h^2\Delta_\gamma)$. The equality $v_h=\psi(h^2\Delta_\gamma)v_h$ means that all frequencies (computed with respect to $\Delta_\gamma$) of $v_h$ are comparable to $-h^{-2}$.

We use a decomposition of $v_h$ as $v_h=v_h^1+v_h^2+v_h^3+v_h^4$ where
\begin{align*}
v_h^1=(1-\chi_0(b_0hD_y))v_h,&\qquad  v_h^2=(\chi_0(b_0hD_y)-\chi_0(b_0^{-1}hD_y))v_h \\
v_h^3=(\chi_0(b_0^{-1}hD_y)-\chi_0(h^{\epsilon}D_y))v_h,& \qquad v_h^4=\chi_0(h^{\epsilon}D_y)v_h,
\end{align*}
that is, a decomposition according to the dual Fourier variable of $y$. Here, $b_0\ll 1$ and $\eps$ will be fixed later. Choosing a good cut-off $\chi_0\in C_c^\infty(\R)$, the proof consists in showing that $\|v_h^j\|_{L^2(M)}=o(1)$ for $j=1,2,3,4$, which contradicts \eqref{e:seekcontrad}. The methods used for each $j$ are quite different, and roughly correspond to the different behaviours of geodesics according to their momentum $\eta\sim D_y$. More precisely:
\begin{itemize}
\item $v^1_h$ corresponds to large $|D_y|$ (i.e., large $\eta$-momenta). For example, the initial data of the vertical Gaussian beams constructed in Section \ref{s:vgb} satisfy $v_h=v_h^1$. To prove that $\|v_h^1\|_{L^2(M)}=o(1)$, we use the positive commutator method. This method dates back at least to \cite[Section 3.5]{hormander1971existence} and has been widely used, for example for proving propagation of singularities for the wave equation. Here, it is based on the relation
\begin{equation*}
[\Delta_\gamma,x\partial_x+(\gamma+1)y\partial_y]=2\Delta_\gamma.
\end{equation*}
Using this relation and computed in two different ways the expression
\begin{equation*}
([h^2\Delta_\gamma+1,\chi_{\epsilon}(x)\phi(y)(x\partial_x+(\gamma+1)y\partial_y)]v_h^1,v_h^1)_{L^2(M)}
\end{equation*}
for some well-chosen cut-offs $\chi_{\epsilon}$ and $\phi$, it is possible to deduce that $\|v_h^1\|_{L^2(M)}=o(1)$.
\item $v^2_h$, $v^3_h$ and $v_h^4$ are ``microlocalized'' in the elliptic part of the symbol of $-\Delta_\gamma$. In some sense, this implies that if we restrict to functions $v_h$ such that $v_h^1=0$, a better resolvent estimate should hold, showing observability of the Schr\"odinger equation in any time, as in the Riemannian case.
\item To prove that $\|v^2_h\|_{L^2(M)}=o(1)$, we consider a defect measure associated to $(v_2^h)_{h>0}$. The key point is that it is invariant along geodesics of the sub-Riemannian metric which reach $\omega$ in finite time (it gives no mass to other geodesics), thus it is null since $\|v^2_h\|_{L^2(\omega)}=o(1)$. This argument is similar in spirit to the construction done in \cite{lebeau1992controle} for the Riemannian Schr\"odinger equation: $v_2^h$ corresponds to the truly elliptic (or semi-classical) regime.
\item The defect measures associated to $v_h^3$ and $v_h^4$ are more complicated to handle since they are invariant along the horizontal geodesics $y=\text{const}$ (which do not enter $\omega$). In some sense, being an approximate solution of $(h^2\Delta_\gamma+1)v=0$ excludes the possibility of ``scarring'' along a single horizontal geodesic; equivalently, for such quasimodes, the mass of $v_h^3$ and $v_h^4$ in $M$ is controlled by their mass in $\omega$. To establish such properties, several tools are available: for example, the papers \cite{burq2012control} and \cite{anantharaman2014semiclassical} deal with similar issues for the elliptic Schr\"odinger equation in tori,  using either normal form arguments or $2$-microlocal techniques. The argument we present in  \cite{letrouit2020observability} relies on the positive commutator method and a similar normal form argument as in  \cite{burq2012control}.
\end{itemize}

\subsection{Further comments and open questions} Together with Corollary 2 in \cite{duyckaerts2012resolvent}, the resolvent estimate proved in Theorem \ref{t:resgamma} implies an observability result for heat-type equations:
\begin{cor} \label{c:corheat} Assume that $\gamma\geq 1$ and that $\omega$ contains a horizontal strip $(-1,1)_x\times I$. For any $s>\frac12(\gamma+1)$ and any $T_0>0$, observability for the heat equation with Dirichlet boundary conditions
\begin{equation} \label{e:heatgam}
\left\lbrace \begin{array}{l}
\partial_tu+(-\Delta_\gamma)^s u=0 \\
u_{|t=0}=u_0 \in L^2(M) \\
u_{|x=\pm 1}=0
\end{array}\right.
\end{equation}
holds in time $T_0$. In other words, there exists $C>0$ such that for any $u_0\in L^2(M)$, there holds
\begin{equation*}
\|e^{-T_0(-\Delta_\gamma)^s}u_0\|_{L^2(M)}^2\leq \int_0^{T_0} \|e^{-t(-\Delta_\gamma)^s}u_0\|_{L^2(\omega)}^2dt.
\end{equation*}
\end{cor}

Let us finally mention a few open questions raised by our study:
\begin{itemize}
\item What happens if $0<\gamma<1$, a case which is not covered by Theorem \ref{t:main2}? Indeed, the sub-Laplacian $\Delta_\gamma$ on $\R\times\T$ is not essentially self-adjoint (see \cite{boscain2016Self}), which means that the Schr\"odinger evolution is not well defined if we do not impose any additional boundary condition on $\{x=0\}$.
\item Is it possible to generalize Theorem \ref{t:main2} to other sub-Laplacians? It seems so that Point (3) might be the easiest one to generalize: very roughly, as seen in Sections \ref{s:seul} and \ref{s:vgb}, some kind of ``normal form'' or ``approximation'' argument for the sub-Laplacian could be relevant since $s<\frac12(\gamma+1)$ leaves some space for perturbative arguments.
\item Is it possible to construct in a more robust way solutions of subelliptic  Schr\"odinger-type equations which are microlocally concentrated in the cone where the principal symbol of the sub-Laplacian vanishes? 
\end{itemize}

\section{Subelliptic Schr\"odinger equation via non-commutative harmonic analysis} \label{s:fermanian}
The last result we present in this survey, namely Theorem \ref{t:main3}, may seem difficult to approach because of the massive use of noncommutative Fourier analysis all along the statement and proofs. Our goal here is to give some keys which could facilitate reading, and to explain why this theory is adapted to analyze subelliptic PDEs. 

The plan of this section is as follows: according to the theory of non-commutative harmonic analysis, we define in Section \ref{s:noncomFou} the (operator-valued) Fourier transform \eqref{e:fouriertransform}, based on the unitary irreducible representations of the group, recalled in \eqref{eq_widehatG}, which form an analog to the usual frequency space. Then, in Section \ref{s:semimeasures}, we follow the same path as for the usual Fourier transform: we use the Fourier inversion formula \eqref{inversionformula} to define in~\eqref{e:defsigma} a class of symbols and the associated semi-classical pseudo-differential operators in~\eqref{def:pseudo}. From this, we deduce the existence of semi-classical measures, which have additional invariance properties when they are associated to solutions of the Schr\"odinger equation. This allows to prove the first part of Theorem \ref{t:main3}. The second part of this theorem is proved using solutions of the Schr\"odinger equation which propagate along the vertical direction: although this is very close in spirit to the vertical Gaussian beam construction of Section \ref{s:vgb}, we have developed in \cite{kammerer2020observability} a more robust framework for these ``wave-packet solutions'', based on the same non-commutative harmonic analysis tools as before.

The material of Sections \ref{s:noncomFou} and \ref{s:semimeasures} borrows many ideas from \cite{fermanian2019quantum} (and of course from \cite{kammerer2020observability}). In the sequel we use the notations of Section \ref{s:thirdmain}.

Additionally, note that the element $(p, q, z) = (p_1, \ldots , p_d, q_1, \ldots , q_d, z)$ of $\H^d$ can be written
\begin{equation*}
(p,q, z)=\exp_{\H^d}(p_1X_1+\ldots+p_dX_d+q_1Y_1+\ldots+q_dY_d+zZ)
\end{equation*}
and with the Baker-Campbell-Hausdorff formula we recover the group law given in Example \ref{exHeis}.

\subsection{Noncommutative Fourier analysis} \label{s:noncomFou} We aim at defining a Fourier transform adapted to $M$, and at proving the associated Fourier inversion formula. This is standard, and the main references are \cite{corwin2004representations} and \cite{taylor1986noncommutative} (see also Appendix A of \cite{kammerer2020observability}). 

The usual Fourier transform $\widehat{f}(\lambda)=\int_{\R^d} f(x)e^{-ix\lambda}dx$ is replaced in this non-commutative setting by the formula
\begin{equation} \label{e:fouriertransform}
{\mathcal F}f(\lambda):=\int_{\H^d} f(x)\left( \pi^{\lambda}_{x }\right)^* \, dx.
\end{equation}
To give it a sense, we need to explain how to define $dx$ and to give a sense to $\lambda$ and $\left( \pi^{\lambda}_{x }\right)^*$.

To define $dx$, we recall that in Section \ref{s:thirdmain}, we defined a scalar product on $\mathfrak z$, for which $\partial_z$ has norm $1$. We also need  a scalar product on $\mathfrak v$, which is obtained by saying that the $2d$ vector fields
\begin{equation*}
X_j=\partial_{x_j}+\frac{y_j}{2}\partial_{z}, \ \  \ Y_j=\partial_{y_j}-\frac{x_j}{2}\partial_{z}, \qquad j=1,\ldots,d
\end{equation*}
form an orthonormal basis. This allows to consider the Lebesgue measure $dv\, dz$ on ${\mathfrak g}={\mathfrak v}\oplus{\mathfrak z}$. Via the identification of $\H^d$ with ${\mathfrak g}$ by the exponential map, this induces a Haar measure $dx$ on $\H^d$ and on its quotient $M$. Thus, we can integrate functions on $\H^d$, define Lebesgue spaces, etc as usual.

Finally, in formula \eqref{e:fouriertransform}, $\left( \pi^{\lambda}_{x }\right)^*$ is an operator. In particular, ${\mathcal F}f(\lambda)$ is operator-valued. The operator $\left( \pi^{\lambda}_{x }\right)^*$ is the adjoint of the  irreducible representations $ \pi^{\lambda}_{x}$ of~$\H^d$. 

\begin{definition}
The (strongly continuous) unitary representations of a locally compact topological group $G$ are the homomorphisms~$\pi:G\rightarrow U(\mathcal{H}_\pi)$ where $U(\mathcal{H}_\pi)$ is the group of unitary operators on a Hilbert space $\mathcal{H}_\pi$, which satisfy that $g\mapsto \pi(g)\xi$ is a norm continuous function for every $\xi\in\mathcal{H}_\pi$. \\
A unitary representation is called irreducible if the only closed linear subspaces of ${\mathcal H}_\pi$ invariant under $\pi(g)$ for all $g\in G$ are $0$ and $\mathcal{H}_\pi$. The set of all unitary irreducible representations (modulo unitary equivalence) is denoted by $\widehat G$.  
\end{definition}
A specificity of Heisenberg groups, and more generally of H-type groups, is that their irreducible representations can be explicitly computed, thanks to Kirillov's theory. We will neither enter the details of Kirillov's theory nor show the computations specific to H-type groups (see Appendix A of \cite{kammerer2020observability}), but only give the explicit expression of irreducible representations:

\begin{itemize}
\item For $\lambda \in \mathfrak z^* \setminus\{0\}\sim \R\setminus \{0\}$ and $x=(p,q,z)\in\H^d$, we consider the operator $\pi_{x}^\lambda$ defined by 
\begin{equation}\label{def:pilambda}
\pi^{\lambda}_{x} \Phi(\xi)=
 {\rm e}^{i(\lambda(z)+\frac12 |\lambda| p\cdot q+\text{sgn}(\lambda)\sqrt{|\lambda|}\xi \cdot q)} \,\Phi \left(\xi+\sqrt{|\lambda|}p\right),
 \end{equation}
which acts on functions $\Phi\in L^2(\R^d)$ (in \eqref{def:pilambda}, $\text{sgn}$ is the sign function).  
Then, $\pi^\lambda(\cdot)=\pi^\lambda_{\cdot}$ defines a unitary representation of $\H^d$ which is infinite dimensional (since the Hilbert space $L^2(\R^d)$ has infinite dimension) and which can be proved to be irreducible.
\item For $\omega\in \mathfrak v ^*$ and $x={\rm Exp} (V+Z)\in \H^d$ with $V\in{\mathfrak v}$ and $Z\in{\mathfrak z}$, we set
\begin{equation} \label{e:finitedim}
\pi^{0,\omega}_x= {\rm e}^{i \omega(V)}
\end{equation}
and $\pi^{0,\omega}(\cdot)=\pi^{0,\omega}_\cdot$ can thus be seen as a $1$-dimensional representation over the Hilbert space $\mathcal{H}_{(0,\omega)}=\mathbb{C}$.
\end{itemize}

 Then, the set $\widehat{\H^d}$ of all unitary irreducible representations modulo unitary equivalence
is parametrized by $({\mathfrak z}^*\setminus \{0\})\sqcup {\mathfrak v}^*$ (see Proposition 2.4 and Theorem 2.5 in Chapter 1 of \cite{taylor1986noncommutative}):
\begin{equation}
\label{eq_widehatG}	
\widehat{\H^d} =
\{\mbox{class of} \ \pi^\lambda \ : \ \lambda \in \mathfrak z^* \setminus\{0\}\}
\sqcup \{\mbox{class of} \ \pi^{0,\omega} \ : \ \omega \in \mathfrak v^* \}.
\end{equation}

The formula \eqref{e:fouriertransform} defines the Fourier transform for $\lambda \in \mathfrak z^* \setminus\{0\}$, but we need to complete it with a formula for $\omega\in \mathfrak v^*$:
 $$
\widehat f(0,\omega)=
 \mathcal F f (0,\omega)
 := \int_{\H^d} f(x) (\pi^{(0,\omega)}_x)^* dx.
$$

With these definitions at hand, one can prove an inversion formula for~$ f \in {\mathcal S}(\H^d)$ and~$x\in \H^d$:
\begin{equation}
\label{inversionformula} f(x)
= c_0 \, \int_{\mathfrak z^* \setminus\{0\}} {\rm{Tr}} \, \Big(\pi^{\lambda}_{x} {\mathcal F}f(\lambda)  \Big)\, |\lambda|^d\,d\lambda \,,
\end{equation}
where ${\rm Tr}$ denotes the trace of operators of ${\mathcal L}(L^2(\R^d))$ (see \cite[Theorem 2.7 in Chapter 1]{taylor1986noncommutative}). Note that in this inversion formula, the finite dimensional irreducible representations given by \eqref{e:finitedim} are absent. 

\subsection{Symbols and semi-classical measures} \label{s:semimeasures}
Starting from the Fourier inversion formula \eqref{inversionformula}, we define symbols (of pseudo-differential operators) on $M$ as a class of functions on $M\times \widehat{\H^d}$.  The set $M\times \widehat{\H^d}$ is interpreted as the phase space of $M$, in analogy with the fact that $\T^d\times \R^d$ is the phase space of the torus $\T^d$.

To start, we note that the set of functions on the quotient $M=\widetilde{\Gamma}\backslash \textbf{H}^d$ is in one-to-one relation with the set of $\widetilde{\Gamma}$-left periodic functions on $\H^d$, i.e., functions on $\H^d$ such that
$$\forall x\in \H^d,\;\;\forall \gamma\in \widetilde{\Gamma} ,\;\; f(\gamma x)=f(x).$$

We consider the class of symbols ${{\mathcal A}_0}$ of fields of operators defined on $M\times \widehat{\H^d}$ by 
$$\sigma(x,\lambda)\in{\mathcal   L}(L^2(\R^d)),\;\;(x,\lambda)\in M\times \widehat{\H^d},$$
that are smooth in the variable $x$ and   Fourier transforms of  functions of the set~${\mathcal S}(\H^d)$ of Schwartz functions on~$\H^d$
 in the variable $\lambda$: for all $(x,\lambda)\in M\times \widehat{\H^d}$,
\begin{equation}\label{e:defsigma}
\sigma(x,\lambda) = {\mathcal F} \kappa_x (\lambda)
\end{equation}
where $\kappa_\cdot(\cdot)\in\mathcal C^\infty(M,\mathcal S(\H^d))$. One can easily check that these symbols form an algebra (which is a motivation for introducing them as Fourier transforms of functions $\kappa_x$).

There is a natural family of dilations on $\H^d$ defined as 
\begin{equation*}
\delta_\varepsilon(x,y,z)=(\varepsilon x,\varepsilon y, \varepsilon^2z), \quad (x,y,z)\in\H^d, \ \varepsilon>0.
\end{equation*}

If $\varepsilon>0$, we associate with $\kappa_x$ (and thus with $\sigma(x,\lambda)$) the function $\kappa_x^\varepsilon$ defined on $\H^d$ by
$$\kappa^\varepsilon_x(\cdot)= \varepsilon^{-Q} \kappa_x(\delta_{\varepsilon^{-1}}(\cdot)),$$
We then  define the semi-classical pseudo-differential operator ${\rm Op}_\varepsilon(\sigma)$ via the identification of functions $f$ on $M$ with $\widetilde{\Gamma}$-left periodic functions on $\H^d$ recalled above:
\begin{equation}\label{def:pseudo}
{\rm Op}_\varepsilon(\sigma)f(x)=  \int_{\H^d} \kappa_x^\varepsilon(y^{-1} x) f(y)dy.
\end{equation}
The operator ${\rm Op}_\varepsilon(\sigma)$ is well-defined as an operator on $M$. Indeed,
\begin{equation*}
{\rm Op}_\varepsilon(\sigma) f(\gamma x)= {\rm Op}_\varepsilon(\sigma) f( x).
\end{equation*}

It is not difficult to check that these definitions yield a ``good symbolic calculus'': for example, the family of operators $\left({\rm Op}_\varepsilon(\sigma)\right)_{\varepsilon >0}$ is uniformly bounded in ${\mathcal L}(L^2(M))$. This allows to define semi-classical measures (see \cite{gerard1991microlocal} and \cite{burq1997mesures} for similar propositions in the Euclidean framework):

\begin{prop}\label{prop:semiclas}
Let $(u^\varepsilon)_{\varepsilon>0}$ be a bounded family in $L^\infty(\R,L^2(M))$. There exist a sequence $(\varepsilon_k)\in(\R_+^*)^{\N}$ with $\varepsilon_k \underset{k\rightarrow +\infty}{\longrightarrow} 0$ and a map $t\mapsto \Gamma_t d\gamma_t$ in $L^\infty(\R,\widetilde{\mathcal{M}}_{ov}^+(M\times \widehat{\H^d}))$ such that we have for all $\theta\in L^1(\R)$ and $\sigma\in\mathcal{A}$,
\begin{equation} \label{e:timeintegratedmeasure}
\int_\R \theta(t)({\rm Op}_{\varepsilon_k}(\sigma)u^{\varepsilon_k}(t),u^{\varepsilon_k}(t))_{L^2(M)} dt  \underset{k\rightarrow +\infty}{\longrightarrow} \int_{\R\times M\times\widehat{\H^d}} \theta(t){\rm Tr}(\sigma(x,\lambda)\Gamma_t(x,\lambda))d\gamma_t(x,\lambda)dt.
\end{equation}
Here $\mathcal{A}$ is the ``closure'' of~$\mathcal A_0$ in some sense made precise in \cite{fermanian2019quantum}.
\end{prop}
Here are a few comments on this statement:
\begin{itemize}
\item We denote by $\widetilde{\mathcal M}_{ov}(M\times \widehat{\H^d})$ the set of pairs $(\gamma,\Gamma)$ where $\gamma$ is a positive Radon measure on~$M\times \widehat{\H^d}$ and $\Gamma=\{\Gamma(x,\lambda)\in {\mathcal L}(L^2(\R^d)):\lambda \in \widehat{\H^d}\}$ is a measurable field of trace-class operators such that
\begin{equation*}
\|\Gamma d \gamma\|_{\mathcal M}:= \int_{M\times \widehat{\H^d}}{\rm Tr}(|\Gamma(x,\lambda)|)d\gamma(x,\lambda)<\infty.
\end{equation*}
Here, as usual, $|\Gamma|:=\sqrt{\Gamma\Gamma^*}$. Finally, we say that
a pair $(\gamma,\Gamma)$ in $ \widetilde {\mathcal M}_{ov}(M\times \widehat{\H^d})$ is {positive} when $\Gamma(x,\lambda)\geq 0$ for $\gamma$-almost all $(x,\lambda)\in M\times \widehat{\H^d}$. In  this case, we write  $(\gamma,\Gamma)\in  \widetilde {\mathcal M}_{ov}^+(M\times \widehat{\H^d})$.
\item This set $\widetilde{\mathcal M}_{ov}(M\times \widehat{\H^d})$ can be identified to the topological dual of the algebra of symbols $\mathcal{A}$ (to be rigorous, it requires to take the quotient by a relation of equivalence on the pairs $(\gamma,\Gamma)$). This is why it naturally appears as a limit of the left-hand side of \eqref{e:timeintegratedmeasure}.
\item Hence, the semi-classical measures that we consider here are operator-valued, whereas semi-classical measures are mostly scalar in the literature, see for example~\cite{lebeau1992controle}. This operator-valued feature is fundamental since it is due to non-commutativity of nilmanifolds, and it is a  consequence of the original features of Fourier analysis on nilpotent groups seen in Section \ref{s:noncomFou}. 
\item The integral in the time variable in \eqref{e:timeintegratedmeasure} may seem weird at first sight, and indeed it is possible to define semi-classical measures of functions on $M$ (thus, time-independent). However, our goal here is to study semi-classical measures associated to solutions of the non-semi-classical Schr\"odinger equation~\eqref{e:Schrod} (a semi-classical Schr\"odinger equation would have the form $i\varepsilon\partial_tu+\varepsilon^2\Delta u=0$). For such equations, it is difficult to derive results for the semi-classical measures at each time $t$ (see also \cite{anantharaman2014semiclassical}). However, one can prove results for the time-averaged semi-classical measures, and this is why we define these time-averaged measures in Proposition \ref{prop:semiclas}.
\item Sections \ref{s:noncomFou} and \ref{s:semimeasures} up to Proposition \ref{prop:semiclas} could be generalized to all graded Lie groups through the  generalization of the tools we use (see Remarks~3.3 and~4.4 in~\cite{fermanian2019semi}). However, the next proposition, namely Proposition \ref{p:measure0}, is specific to H-type groups (in particular, Heisenberg groups) since its proof uses the explicit expressions of irreducible representations seen in Section \ref{s:noncomFou}.
\end{itemize}

The semi-classical measures associated (by Proposition \ref{prop:semiclas}) to families of solutions to the Schr\"odinger equation~\eqref{e:Schrod} have special features, which are the subject of Proposition \ref{p:measure0} below. To state it properly, we need some definitions. 

In the (non compact) group~$\H^d$, the operator 
$$H(\lambda)=|\lambda| \sum_{j=1}^{d} \left( -\partial_{\xi_j}^2+\xi_j^2\right)$$ is the Fourier resolution of the sub-Laplacian $-\Delta_{\H^d}$ above $\lambda\in \mathfrak z^*\setminus \{0\}$, meaning that 
$$\forall f\in{\mathcal S}(\H^d),\;\;\mathcal F(-\Delta_{\H^d} f)(\lambda)= H(\lambda) {\mathcal F} (f) (\lambda).$$
Up to a constant, this is a quantum harmonic oscillator with 
discrete spectrum $\{|\lambda|(2n+d), n\in \N\}$ 
and finite dimensional eigenspaces.
For each eigenvalue $|\lambda|(2n+d)$, we denote by $\Pi_n^{(\lambda)}$ and $\mathcal V^{(\lambda)}_n$ the corresponding spectral orthogonal projection and eigenspace.
 
\begin{prop} \label{p:measure0}
Assume $\Gamma_t d\gamma_t$ is associated with a family of solutions to~\eqref{e:Schrod}.
\begin{enumerate}
\item For $(x,\lambda)\in  M\times \mathfrak z^*$
\begin{equation}\label{eq:decomp}
\Gamma_t(x,\lambda)=\sum_{n\in\N } \Gamma_{n,t}(x,\lambda)\;\;{ with}\;\; \Gamma_{n,t}(x,\lambda):= \Pi_n^{(\lambda)}\Gamma_t(x,\lambda) \Pi_n^{(\lambda)}.
\end{equation}
Moreover,
the map $(t,x,\lambda)\mapsto \Gamma_{n,t}(x,\lambda) d\gamma_t(x,\lambda)$ defines
 a continuous function from  $\R$ into the set of distributions on $M\times (\mathfrak z^*\setminus\{0\})$ valued in the finite dimensional space ${\mathcal L}(\mathcal V_n^{(\lambda)})$ which
satisfies
\begin{equation}\label{transport}
\left(\partial_t -(n+\frac d 2) {\mathcal Z}^{(\lambda)}\right)\left(\Gamma_{n,t}(x,\lambda)d\gamma_t(x,\lambda)\right)=0
\end{equation}
\item For $(x,(0,\omega))\in M\times \mathfrak v^*$, the scalar measure $\Gamma_t d\gamma_t$ is invariant under the flow
$$\Xi ^s:(x,\omega) \mapsto (x{\rm Exp} (s\omega\cdot V),\omega).$$
Here, $\omega\cdot V=\sum_{j=1}^{2d} \omega_jV_j$ where $\omega_j$ denote the coordinates of $\omega$ in the dual basis of $V$.
\end{enumerate}

\end{prop}

The relations \eqref{eq:decomp} and \eqref{transport} describe the same phenomenon as \eqref{e:decompoL2grushin}: the Schr\"odinger equation behaves as a superposition of waves traveling at different speeds along the vertical axis.

\subsection{Proof of Theorem \ref{t:main3}} We finally explain how Theorem \ref{t:main3} follows from Proposition \ref{p:measure0}.
\subsubsection{Proof of Point (1) of Theorem \ref{t:main3}} As in Section \ref{s:resgrushin}, to prove Point (1) of Theorem \ref{t:main3}, it is sufficient to prove a \emph{spectrally localized} observability inequality.
Let $h>0$ and $\psi\in C_c^\infty((1/2,2),[0,1])$. Using functional calculus, we set
 \begin{equation}\label{def:Pih}
 \mathcal{P}_h f= \psi\left(-h^2\left(\frac12 \Delta_M+\mathbb{V}\right)\right) f,\;\; f\in L^2(M).
 \end{equation}
 We seek to prove
 \begin{equation}\label{obs_loc}
\| \mathcal{P}_h u_0\|^2_{L^2(M)} \leq C_0 \int_0^T \left\|   {\rm e}^{it(\frac12 \Delta_M+\mathbb{V})}  \mathcal{P}_h u_0\right\|^2_{L^2(U)} dt.
\end{equation}

We argue by contradiction. If~\eqref{obs_loc} is false, then there exist  $(u_0^k)_{k\in\N}$ and $(h_k)_{k\in\N}$ such that $u_0^k= \mathcal{P}_{h_k} u_0^k$,
\begin{equation}\label{hyp}
\| u^k_0\|_{L^2(M)}=1\;\;\mbox{and}\;\;  \int_0^T \left\|   u_k(t)\right\|_{L^2(U)}^2 dt\underset{k\rightarrow +\infty}{\longrightarrow}0.
\end{equation}
where
$$u_k(t)=   {\rm e}^{it(\frac12 \Delta_M+\mathbb{V})} \mathcal{P}_{h_k} u_0^k.$$
We consider (after extraction of a subsequence if necessary), the semi-classical measure $\Gamma_t d\gamma_t$ of $u_k(t)$ given by Proposition~ \ref{prop:semiclas} and satisfying the properties listed in Proposition~\ref{p:measure0}. The goal is to prove that $\gamma_t\equiv 0$.  

Using the second part of \eqref{hyp}, one obtains
\begin{equation*} 
\int_0^T \int_{U\times \widehat{\H^d}} {\rm Tr} (\Gamma_t(x,\lambda)) d\gamma_t (x,\lambda)dt=0
\end{equation*}
i.e., $\gamma_t\equiv 0$ above $U$. Setting $\gamma_{n,t} (x,\lambda)= {\rm Tr}\left(\Gamma_{n,t}(x,\lambda)\right) \gamma_t(x,\lambda)$, and using the positivity of $\Gamma_t$, one can deduce that
\begin{equation}\label{e:nullityoverU} 
\int_{U\times \mathfrak z^*} d\gamma_{n,t} (x,\lambda)=0,\;\; \text{for almost every } t\in [0,T],\;\; \forall n\in\N.
\end{equation}
The transport equation~\eqref{transport} tells us that $\gamma_{n,t}$ travels at speed $n+\frac{d}{2}$ along the $z$-axis, hence not slower than $\gamma_{0,t}$. Using (H-GCC) together with \eqref{e:nullityoverU}, we get that $\gamma_{n,t}\equiv 0$ for any $n\in\N$, hence $\gamma_t\equiv 0$. This contradicts the conservation of energy (i.e., the first part of \eqref{hyp}). Thus, \eqref{obs_loc} is proved.

\begin{remark} \label{r:comparsuncfk}
It is tempting to compare the approaches developed in Sections \ref{s:sun} and \ref{s:fermanian}, which share the common goal of proving observability inequalities for subelliptic Schr\"odinger equations. On one side, the semi-classical measures of Section \ref{s:fermanian} seem particularly adapted: once defined the operator-valued Fourier transform, the definitions of symbols and semi-classical measures are natural since they are modeled on the Euclidean case. But this approach has the drawback to require the knowledge of global objects on the manifold (the representations), and for the moment their local and geometric aspects are not sufficiently well understood to handle more general geometric situations.
On the contrary, the proof of Section \ref{s:sun}, which uses usual pseudo-differential tools (i.e., Euclidean or Riemannian ones), is directly linked with the underlying geometry (see Section \ref{s:proofres}) but we see the limits of these tools in the fact that already for the simple Baouendi-Grushin models, the computations are sophisticated and sometimes heavy.
\end{remark}
\subsubsection{Proof of Point (2) of Theorem \ref{t:main3}} As for Point (3) of Theorem \ref{t:main2}, to disprove the observability inequality \eqref{obs}, we construct a particular family of solutions of the Schr\"odinger equation \eqref{e:Schrod}, called wave packets (and which are somehow related to the vertical Gaussian beams of Section \ref{s:vgb}).

In the Euclidean context, given $(x_0,\xi_0)\in\R^d\times\R^d$ and $a\in\mathcal S(\R^d)$, the associated wave packet is the family (indexed by $\varepsilon$) of functions
\begin{equation}\label{WPeuclid}
u_{\rm eucl}^\varepsilon(x)= \varepsilon^{-d/4} a\left(\frac{x-x_0}{\sqrt\varepsilon}\right) {\rm e}^{\frac i\varepsilon \xi_0\cdot(x-x_0)}, \;\; x\in\R^d.
\end{equation}
The oscillation along $\xi_0$ is forced by the term ${\rm e}^{\frac i\varepsilon \xi_0\cdot(x-x_0)}$ and the concentration on $x_0$ is performed via $a(\cdot/\sqrt\varepsilon)$ at the scale $\sqrt\varepsilon$ for symmetry reasons:  the $\varepsilon$-Fourier transform of $u^\varepsilon_{\rm eucl}$,
$\varepsilon^{-d/2}\widehat u^\varepsilon_{\rm eucl} (\xi/\varepsilon)$ presents a concentration on $\xi_0$ at the scale $\sqrt\varepsilon$. Taking $a$ compactly supported in the interior of a unit cell  for the torus, one can generalize their definition to the case of the torus by extending them by periodicity.  

To perform a similar construction in the non-commutative setting, more precisely in $\H^d$, we replace $a(\cdot/\sqrt\varepsilon)$ by
$$a_\varepsilon(x)=a \left(\delta_{\varepsilon^{-1/2 }}(x )\right)$$
for some $a\in\mathcal C_c^\infty(G)$, and the oscillations ${\rm e}^{\frac i\varepsilon \xi_0\cdot(x-x_0)}$ by
$$e_\varepsilon(x)=\left(\pi^{\lambda_\varepsilon}_x \Phi_1, \Phi_2\right), \;\;\lambda_\varepsilon= \frac{\lambda_0}{\varepsilon^2}$$
where $\lambda_0\in\mathfrak z^*$ and $\Phi_1,\; \Phi_2\in{\mathcal S}(\R^d)$.

Using the multiplication on the left by elements of $\widetilde\Gamma$, one can define a periodization operator $\mathbb P$ which associates to functions on $\H^d$ whose support is contained in a unit cell of $M$ the $\widetilde{\Gamma}$-left periodic function obtained by periodization. We restrict to $\varepsilon\in (0,1)$ We consider a unit cell of $M$, i.e., a subset~$\mathcal B$ of~$\H^d$ which is a neighborhood of~$1_{\H^d}$ and such that~$\cup_{\gamma\in \widetilde{\Gamma}} (\gamma {\mathcal B} )=\H^d$ and we choose functions~$a$ that are in~${\mathcal C}^\infty_c(\mathcal B)$.

\begin{prop}\label{prop:WP}
Let $\Phi_1,\Phi_2\in \mathcal S(\R^d)$, $a\in C_c^\infty(\mathcal B)$, 
$x_0\in M$,  $\lambda_0\in\mathfrak z^*\setminus\{0\}$. Then, there exists $\varepsilon_0>0$ such that  the family~$(v^\varepsilon)_{\varepsilon\in(0,\varepsilon_0)}$  defined by
$$v^\varepsilon (x)=|\lambda_\varepsilon|^{d/2}\,\varepsilon^{-p/2}
\, {\mathbb P} (e_\varepsilon a_\varepsilon)( x_0^{-1} x),$$
has only one  semi-classical measure $\Gamma d\gamma$  where
\begin{equation}\label{eq:WPmeasures}
 \gamma=c_a\, \delta(x-x_0) \otimes \delta (\lambda-\lambda_0),\;\; c_a= \|\Phi_2\|^2\int_{G_{\mathfrak z}} |a( x_{\mathfrak z})|^2 d x_{\mathfrak z},
 \end{equation}
 and $\Gamma$ is the operator
defined by
$$\Gamma\Phi= \frac{(\Phi,\Phi_1)}{\|\Phi_1\|^2} \Phi_1,\;\;\forall \Phi\in L^2(\R^d).$$
\end{prop}

What is really important in the above proposition is that $\gamma$ is concentrated on a single point of the phase space.

In the rest of the proof, we say that the family $v^\varepsilon$ is a wave packet on $M$ with cores $(x_0,\lambda_0)$, profile~$a$ and harmonics $(\Phi_1,\Phi_2)$, and we write
$$v^\varepsilon= WP^\varepsilon_{x_0,\lambda_0}(a,\Phi_1,\Phi_2)=|\lambda_\varepsilon|^{d/2}\,\varepsilon^{-p/2}
\, {\mathbb P} (e_\varepsilon a_\varepsilon)(  x_0^{-1}x).$$
\begin{remark}
In \cite{kammerer2020observability} (see Section 4 and Appendix C), we develop a more general theory of wave packets, notably showing that the structure of wave packets is preserved by the evolution under the Schr\"odinger flow.
\end{remark}

We take as initial data in~\eqref{e:Schrod} a wave packet $u^\varepsilon_0$  in~$M$ with harmonics given by the first Hermite function $h_0$:
$$u^\varepsilon_0=WP^\varepsilon_{x_0,\lambda_0}(a,h_0,h_0).$$
We denote by $u^\varepsilon(t)$ the associated solution, $u^\varepsilon(t)=  {\rm e}^{it(\frac12 \Delta_M+\mathbb{V})} u^\varepsilon_0$. Our choice of harmonics for $u^\varepsilon_0$ guarantees that the semi-classical measure $\Gamma_td\gamma_t$ associated to these solutions, when decomposed according to Proposition \ref{p:measure0}, has only one non-vanishing component, which corresponds to $n=0$. In other words,
\begin{equation}\label{eq:gammat}
\gamma_t (x,\lambda)= c \, \delta\left( x-{\rm Exp} \left(t\frac d2\mathcal Z^{(\lambda)} \right)x_0\right)\otimes \delta(\lambda-\lambda_0)
\end{equation}
and $\Gamma_{0,t}$ is the orthogonal projector on $h_0$.

Now, using the assumptions made in Point  (2) of Theorem \ref{t:main3}, there exists a continuous function $\phi:M\rightarrow [0,1]$ such that $\phi( \Phi_0^s(x_0,\lambda_0))=0$ for any $s\in [0,T]$ and $\phi=1$ on $\overline{U}\times \mathfrak z^*$. From this, we deduce
\begin{align*}
0\leq \int_0^{T} \int_U | u^\varepsilon(t,x)|^2 dx dt \leq \int_0^{T} \int_M \phi(x)| u^\varepsilon(t,x)|^2 dx dt \underset{\varepsilon\rightarrow 0}{\longrightarrow} & \int_0^{T} \int_{M\times \mathfrak z^*} \phi(x)d\gamma_t (x,\lambda) dt=0.
\end{align*}
Therefore, the observability inequality \eqref{obs} cannot hold.

\section{Perspectives and open problems} \label{s:open}
In the field of subelliptic PDEs, and notably concerning controllability/observability, several interesting questions remain unanswered, and this concluding section lists a few of them.

\subsection{Observability of the heat equation for other sub-Laplacians.} \label{s:conjheat} As explained in Section \ref{s:subheatequations}, the observability properties of subelliptic heat equations are known only in particular geometries. More general results would require a deeper understanding of the geometric meaning of the solutions constructed in  \cite{beauchard2014null} or \cite{koenig2017non}. Let us formulate two conjectures: 
\begin{enumerate}
\item For any sub-Laplacian of step $2$, if $M\setminus \omega$ has non-empty interior, observability of the associated heat equation fails for sufficiently small times $T>0$;
\item For any sub-Laplacian of step $\geq 3$, if $M\setminus \omega$ has non-empty interior, observability of the associated heat equation fails for any time $T>0$.
\end{enumerate}
These conjectures are inspired by the results mentioned in Section \ref{s:subheatequations} and by the paper \cite{laurent2020tunneling} (see notably Section 1.4).

\subsection{Observability of Schr\"odinger for other sub-Laplacians.} \label{s:openobs} Even in the Riemannian case, the observability properties of the Schr\"odinger equation remain mysterious: although (GCC) is known to be a sufficient condition for observability, it is not a necessary condition (see Section \ref{s:obsSchrodRiem}). In the sub-Riemannian case,  the problem is even ``more open'', since no general sufficient condition is known for the moment, except trivial ones: only very particular geometries have been explored (see Theorems \ref{t:burqsun}, \ref{t:main2} and \ref{t:main3}), and they rely on tools which are not robust enough to cover general (in particular non-flat) sub-Riemannian geometries.

\subsection{Propagation of singularities for subelliptic PDEs.} \label{s:opensing} Observability properties of the wave and Schr\"odinger equations are related to propagation of singularities of their solutions (see for example \cite{bardos1992sharp}). The propagation of singularities for sub-Riemannian wave equations has been addressed in a series of papers by Melrose, R. Lascar and B. Lascar, culminating with the general result of \cite{melrose1984propagation}, which we revisited in \cite{letrouit2021propagation} and \cite{colin2021propagation}.

Besides, the understanding of propagation of singularities for sub-Riemannian Schr\"odinger equations surely requires the introduction of a notion of singularity ``adapted to the sub-Riemannian geometry'', i.e., taking into account the number of brackets needed to generate each direction.

\subsection{Trace formulas} \label{s:opentrace} For most (sub-)Laplacians $\Delta$, we do not have access to the knowledge of the full spectrum, i.e., to all eigenvalues. But it is sometimes possible to compute quantities of the form 
\begin{equation*}
\sum_{n\in\N} f(\lambda_n)
\end{equation*}
where $f$ is a (possibly complex-valued) function and $\lambda_n$ describes the spectrum (with multiplicities) of $-\Delta$, i.e., $-\Delta\varphi_n=\lambda_n\varphi_n$ for smooth functions $\varphi_n$. Classical choices for $f$ are the following: $f(x)=e^{-tx}$ (heat equation), $f(x)=|x|^{-s}$ (zeta functions), $f(x)=\cos(t\sqrt{x})$ (wave equation), $f(x)=e^{-itx/h}$ (semi-classical Schr\"odinger equation). 

The literature on trace formulas in Riemannian manifolds is vast. But in the sub-Rzhak iemannian case, only few trace formulas have been established, and most of them are formulated with the heat kernel. It would be of interest to prove trace formulas for other kernels.

\subsection{Eigenfunctions and quasimodes of sub-Laplacians.}\label{s:openeigen} The properties of eigenfunctions and quasimodes of sub-Laplacians remain widely unknown. Beside the concentration results given by observability properties (see for example Corollary \ref{c:eigenfunctions}), one could expect to characterize the weak limits of high-frequency eigenfunctions (or of the square of their modulus) in the limit where the eigenvalue tends to $+\infty$: these weak limits are known as Quantum Limits and they were widely studied in the Riemannian case.  In \cite{de2018spectral}, the authors undertook their study in the sub-Riemannian case, proving that they concentrate (except for a null-density subsequence) on the characteristic cone $\Sigma$, and also showing a Quantum Ergodicity result valid for 3D contact sub-Laplacians with ergodic Reeb flow. Their study was pursued in \cite{letrouit2020quantum}, in the case where a commutativity assumption on the vector fields involved in the definition of the sub-Laplacian is satisfied: then, techniques coming from joint spectral calculus can be applied. In particular, we characterized all Quantum Limits of a family of sub-Laplacians obtained as a product of flat 3D Heisenberg sub-Laplacians. But, for example, the higher-dimensional (non-flat) contact case remains open.


\bibliographystyle{alpha}
\bibliography{biblioactestsgv2}

\end{document}